\documentclass[twoside,12pt,leqno]{article}
-\headheight\advance\topmargin by -\headsep

\usepackage{amsmath}
\usepackage{amsthm}
\usepackage{amssymb}
\usepackage{amscd}
\usepackage{graphicx}
\usepackage{epsfig}

\numberwithin{equation}{section} \allowdisplaybreaks

\newtheorem{theorem}{Theorem}[section]
\newtheorem{lemma}{Lemma}[section]
\newtheorem{proposition}{Proposition}[section]
\newtheorem{corollary}{Corollary}[section]

\theoremstyle{definition}
\newtheorem{definition}{Definition}[section]
\newtheorem{example}{Example} [section]

\newtheorem{remark}{Remark}[section]

\begin{document}

\newcommand{\TT}{\mathbb{T}}
\newcommand{\ZZ}{\mathbb{Z}}
\newcommand{\CC}{\mathbb{C}}
\newcommand{\RR}{\mathbb{R}}
\newcommand{\cE}{\mathcal{E}}
\newcommand{\cL}{\mathcal{L}}
\newcommand{\cF}{\mathcal{F}}
\newcommand{\bL}{\bar{L}}
\newcommand{\para}{\kern0.2em{\backslash} \kern-0.7em {\backslash} \kern0.2em }



\setcounter{page}{1} \thispagestyle{empty}



\markboth{A.~Weinstein and M.~Zambon}{Variations on
Prequantization}

\label{firstpage} $ $
\bigskip

\bigskip

\centerline{{\Large  Variations on Prequantization}}

\bigskip
\bigskip
\centerline{{\large by  Alan Weinstein and Marco Zambon}}

\vspace*{.7cm}

\begin{abstract}
We extend known prequantization procedures for Poisson and
presymplectic manifolds by defining the prequantization of a Dirac
manifold $P$ as a principal $U(1)$-bundle $Q$ with a compatible
Dirac-Jacobi structure.  We study the action of Poisson algebras
of admissible functions on $P$ on various spaces of locally (with
respect to $P$) defined functions on $Q$, via hamiltonian vector
fields. Finally, guided by examples arising in complex analysis
and contact geometry, we propose an extension of the notion of
prequantization in which the action of $U(1)$ on $Q$ is permitted
to have some fixed points.
\end{abstract}

\pagestyle{myheadings}

\begin{center}
\emph{Dedicated to the memory of Professor Shiing-Shen Chern}
\end {center}

\section{Introduction}

Prequantization in symplectic geometry attaches to a symplectic
manifold $P$ a hermitian line bundle $K$ (or the corresponding
principal $U(1)$-bundle $Q$), with a connection whose curvature
form is the symplectic structure.  The Poisson Lie algebra
$C^\infty(P)$ then acts faithfully on the space $\Gamma(K)$ of
sections of $K$ (or antiequivariant functions on $Q$). Imposing a
polarization $\Pi$ cuts down $\Gamma(K)$ to a smaller, more
``physically appropriate'' space $\Gamma_\Pi(K)$ on which a
subalgebra of $C^\infty(P)$ may still act.  By polarizing and
looking at the ``ladder'' of sections of tensor powers $K^{\otimes
n}$ (or functions on $Q$ transforming according to all the
negative tensor powers of the standard representation of $U(1)$),
one gets an ``asymptotic representation'' of the full algebra
$C^\infty(P)$.  All of this often
 goes under the name of  geometric quantization, with the last
 step closely related to deformation quantization.

For systems with constraints or systems with symmetry,
the phase space $P$ may be a presymplectic or Poisson
manifold.  Prequantization, and sometimes the full procedure of
geometric quantization, has been carried out in these settings by
several authors; their work is cited below.

The principal aim of this paper is to suggest two extensions of
the prequantization construction which originally arose in an
example coming from contact geometry.  The first, which unifies
the presymplectic and Poisson cases and thus permits the simultaneous
application of constraints and symmetry, is to allow $P$ to be a Dirac
manifold.  The second is to allow the $U(1)$ action on $Q$ to have
fixed points when $P$ has a boundary, so that the antiequivariant
functions become sections of a sheaf rather than a line bundle over
$P$.
In the course of the paper, we also make some new observations
concerning the Poisson and presymplectic cases.

\subsection{Symplectic prequantization}
On a symplectic manifold $(P, \omega)$, one defines the
hamiltonian vector field $X_f$ of the function $f$ by
$\omega(X_f,\cdot)=df$, and one has the Lie algebra bracket
$\{f,g\}=\omega(X_f,X_g)$ on $C^{\infty}(P)$.  A closed 2-form
$\omega$ is called integral if its de Rham cohomology class
$[\omega]\in H^2(M,\RR)$ is integral, i.e. if it is in the image
of the homomorphism $i_*:H^2(M,\ZZ)\to H^2(M,\RR)$
associated with the inclusion $i:\ZZ \to \RR$ of coefficient
groups.

When $\omega$ is integral, following Kostant \cite{Ko1}, we
prequantize $(P,\omega)$ by choosing a hermitian line bundle
 $K$
bundle over $P$ with first Chern class  in $i_{*}^{{-1}}[\omega]$.
  Then there is a
 connection $\nabla$ on $K$ with curvature $2 \pi i\omega$.
Associating to each function $f$ the operator $\hat{f}$ on
$\Gamma(K)$ defined by\footnote{Our convention for the Poisson bracket
differs by a sign from that of \cite{Gu} and \cite{Ko1};
consequently our formula for $\hat{f}$ and Equation \eqref{sympl}
below differ by a sign too. Our sign has the property that the map
from functions to their hamiltonian vector fields is an
antihomomorphism from Poisson brackets to Lie brackets.}
 $\hat{f}(s)=-[\nabla_{X_f}s+2\pi i
fs]$, we obtain a faithful Lie algebra representation of
$C^{\infty}(P)$ on $\Gamma(K)$.

The construction above is equivalent to the following, due to
Souriau \cite{So}: let $Q$ be the principal $U(1)$-bundle
associated to $K$. Denote by $\sigma$ the connection form on $Q$
corresponding to $\nabla$ (so $d\sigma=\pi^*\omega$, where $\pi:Q
\rightarrow P$), and by $E$ the infinitesimal generator of the
$U(1)$ action on $Q$.
 We can
identify the sections of $K$
with functions $\bar{s}:Q \rightarrow \CC$ which are
$U(1)$-antiequivariant (i.e. $\bar{s}(x\cdot t)=\bar{s}(x)\cdot
t^{-1}$ for $x\in Q$, $t\in U(1)$, or equivalently
$E(\bar{s})=-2\pi i \bar{s}$),
 and then the
operator
$\hat{f}$ on $\Gamma(K)$ corresponds to the action of  the
vector field
\begin{eqnarray} \label{sympl} -X_f^H+\pi^*fE, \end{eqnarray} where the
superscript $^H$ denotes the horizontal lift to  $Q$ of a vector
field on $P$. Notice that $\sigma$ is a contact form on $Q$ and
that $X_f^H-\pi^*fE$ is just the hamiltonian vector field of
$\pi^*f$ with respect to this contact form (viewed as  a Jacobi
structure; see Section \ref{secjacdir}).

\subsection{Presymplectic prequantization}
Prequantization of a presymplectic
manifolds
 $(P,\omega)$ for which $\omega$ is integral and of constant
rank\footnote{Unlike many other authors (including some of those
cited here), we will use the work ``presymplectic'' to describe
any manifold endowed with a closed 2-form, even if the form does
not have constant rank.} was introduced by G\"unther \cite{Gu} (see also Gotay and
Sniatycki \cite{GoS} and Vaisman \cite{Va1}). G\"unther represents
the Lie algebra of functions constant along the leaves of $\ker
\omega$ by assigning to each such function $f$ the equivalence
class of vector fields on $Q$ given by formula \eqref{sympl},
where $X_f$ now stands for the equivalence class of vector fields
satisfying $\omega(X_f,\cdot)=df$.

\subsection{Poisson prequantization}
Prequantization of Poisson manifolds $(P,\Lambda)$ was first
investigated algebraically  by Huebschmann \cite{Hu}, in terms of
line bundles by Vaisman \cite{Va2},  and then in terms of circle
bundles by Chinea, Marrero, and de Leon \cite{CMdL}. When the
Poisson cohomology class $[\Lambda]\in H^2_{\Lambda}(P)$ is
the image of an integral de Rham class $[\Omega]$  under the map
given by contraction with $\Lambda$, a $U(1)$-bundle $Q$ with
first Chern class in $i_{*}^{-1}[\Omega]$
may be given a Jacobi structure for which
 the map that assigns to $f\in C^{\infty}(P)$
the hamiltonian vector field (with respect to the Jacobi
structure) of $-\pi^*f$ is a Lie algebra homomorphism. This gives a
(not always faithful) representation of $C^{\infty}(P)$.

\subsection{Dirac prequantization}
We will unite the results in the previous two paragraphs by using
 Dirac manifolds.  These
were introduced by Courant \cite{Co}
 and include both Poisson  and presymplectic manifolds as
special cases.   On the other hand, Jacobi manifolds had already
been introduced by Kirillov \cite{Ki} and Lichnerowicz \cite{Li},
including Poisson, conformally symplectic, and contact manifolds
as special cases.  All of these generalizations of Poisson
structures were encompassed in the definition by Wade \cite{Wa} of
Dirac-Jacobi\footnote{Wade actually calls them $\cE^1(M)$-Dirac
manifolds; we will stick to the terminology ``Dirac-Jacobi'', as
introduced in \cite{GrM}.}
  manifolds.

 To prequantize a Dirac manifold $P$, we will impose an
integrality condition on $P$ which implies the existence of a
$U(1)$-bundle $\pi:Q \rightarrow P$ with a connection which will
be used to construct a Dirac-Jacobi structure on $Q$.
Prequantization of  (suitable) functions $g\in C^{\infty}(P)$ is
achieved ``Kostant-style'' by associating to $g$ the equivalence
class of the hamiltonian vector fields of $-\pi^*g$ and by letting
this equivalence class act on a suitable subset of the
$U(1)$-antiequivariant functions on $Q$, or equivalently by
letting $\pi^*g$ act by the bracket of functions on $Q$. The same
prequantization representation can be realized as an action on
sections of a hermitian line bundle over $P$ with an
$L$-connection, where $L$ is the Lie algebroid given by the Dirac
manifold $P$.

We also look at the following very natural example,
discovered by Claude LeBrun. Given a contact manifold $M$ with
contact distribution $C\subset TM$, the nonzero part of its
annihilator $C^{\circ}$ is a symplectic submanifold of $T^*M$.
When the contact structure is cooriented, we may choose the
positive half $C^{\circ}_+$ of this submanifold.
By adjoining to
$C^{\circ}_+$
the ``the section at
infinity of $T^*M$'' we obtain a manifold with boundary, on which
the symplectic structure on  $C^{\circ}_+$ extends to give a Poisson structure.
We call this a ``LeBrun-Poisson manifold''.
If now we additionally adjoin  the zero section of $T^*M$
we obtain a Dirac manifold $P$.

First we will describe the
prequantization  $U(1)$-bundle of $P$,
then
 we will modify it by
collapsing to points the fibers over one of the two  boundary components
 and by applying a conformal
change. At the end, restricting this construction to
the LeBrun-Poisson manifold (which sits as an open set inside $P$),
 we will obtain a contact manifold in which  $M$
 sits as a contact submanifold.

\subsection{Organization of the paper}
In Sections \ref{secdir} and \ref{secjacdir} we collect known
facts about Dirac and Dirac-Jacobi manifolds. In Section
\ref{spaces} we state our prequantization condition and describe
the Dirac-Jacobi structure on the prequantization space of a Dirac
manifold. In Section \ref{representation} we study the
corresponding prequantization representation, and in Section
\ref{line} we derive the same representation by considering
hermitian line bundles endowed with $L$-connections. In Section
\ref{seclebrun} we study the prequantization of LeBrun's examples,
and in Section \ref{fixed} we allow prequantization $U(1)$-bundles
to have fixed points, and we endow them with contact structures.
We conclude with some remarks in Section \ref{remarks}.

\noindent{\textbf{Acknowledgements:}}  A.W. would like to thank the
 Institut Math\'ematique de Jussieu and \'Ecole Polytechnique for
 hospitality while this
  paper was being prepared.   M.Z. is grateful to Xiang
Tang for helpful discussions and advice in the early stages of
this work.  We would both like to thank the local organizers of
 Poisson 2004 for encouraging the writing of this article by insisting
 on the publication of a volume of proceedings (as well as for the
 superb organization of the meeting itself).

\section{Dirac manifolds}\label{secdir}

We start by recalling some facts from \cite{Co}.
\begin{definition}[\cite{Co}, Def 1.1.1]\label{lindirstr}
A {\bf Dirac structure} on a vector space $V$ is a maximal isotropic
subspace $L \subset V\oplus V^*$ with respect to
the symmetric pairing
\begin{eqnarray}\label{symm}  \langle X_1\oplus\xi_1,X_2\oplus\xi_2 \rangle_+
= \frac{1}{2}(i_{X_2}\xi_1
 +i_{X_1}\xi_2). \end{eqnarray}
\end{definition}
$L$ necessarily has the same dimension as $V$,
and denoting by $\rho_V$ and $\rho_{V^*}$ the projections of
$V\oplus V^*$ onto $V$ and  $V^*$ respectively, we have
 \begin{eqnarray} \label{annih}\rho_V(L)=(L \cap V^*)^{\circ}
 \text{     and     }\rho_{V^*}(L)=(L \cap V)^{\circ} \end{eqnarray}
where the symbol $^{\circ}$ denotes the annihilator. It follows
that $L$ induces (and is equivalent to) a skew bilinear form on
$\rho_V(L)$ or a  bivector on $V/L \cap V$ (\cite{Co}, Prop.
1.1.4). If $(V,L)$ is a Dirac vector space and $i:W \rightarrow V$
a linear map, then one obtains a pullback Dirac structure on $W$
by $\{Y\oplus i^*\xi:iY\oplus \xi\in L\}$; one calls a map between
Dirac vector spaces ``backward Dirac map'' if it pulls back the
Dirac structure of the target vector space to the one on the
source vector space \cite{BuR}. Similarly, given a linear map
 $p: V \rightarrow Z$, one obtains
 a pushforward  Dirac structure on $Z$ by
$\{pX\oplus \xi:X\oplus p^*\xi\in L\}$, and one thus has a notion of
``forward Dirac map'' as well.

On a manifold $M$, a maximal isotropic subbundle $L\subset
TM\oplus T^*M$ is called an {\bf almost Dirac structure} on $M$.
The appropriate integrability condition was discovered by Courant
(\cite{Co}, Def. 2.3.1):
\begin{definition}\label{defdirstr}
A {\bf Dirac structure} on $M$ is an almost Dirac structure $L$ on
$M$ whose space of sections is closed under the Courant
bracket on sections of $TM\oplus T^*M$, which is defined by
\begin{eqnarray} \label{coubra}\;\;\;\;\;\;\;
[X_1\oplus\xi_1,X_2\oplus\xi_2]=
\big([X_1,X_2]\;\oplus\;\cL_{X_1}\xi_2-\cL_{X_2}\xi_1+\frac{1}{2}d(
i_{X_2}\xi_1 -i_{X_1}\xi_2 )\big).\end{eqnarray}
\end{definition}

When an almost Dirac structure $L$ is integrable,
$(L,\rho_{TM}|_L,[\cdot,\cdot])$ is a Lie
algebroid\footnote{Recall that a Lie algebroid is a vector bundle
$A$ over a manifold $M$ together with a Lie bracket
$[\cdot,\cdot]$ on its space of sections and a bundle map $\rho:
A\rightarrow TM$ (the ``anchor'') such that the Leibniz rule
$[s_1,fs_2]=\rho s_1(f)\cdot s_2+f\cdot [s_1,s_2]$ is satisfied
for all sections $s_1,s_2$ of $A$ and functions $f$ on $M$.}
 (\cite{Co}, Thm. 2.3.4). The singular
distribution $\rho_{TM}(L)$ is then integrable in the sense of
Stefan and Sussmann \cite{Su} and gives rise to a singular
foliation of $M$. The Dirac structure induces a closed 2-form
(presymplectic form) on each leaf of this foliation (\cite{Co},
Thm. 2.3.6). The distribution $L\cap V$, called the {\bf
characteristic distribution}, is singular in a different way. Its
annihilator $\rho_{T^*M}(L)$ is closed in the cotangent bundle,
but the distribution itself is not closed
unless it has
constant rank. It is not always integrable, either.
 (See Example \ref{notadm} and the beginning of Section
\ref{seclebrun}.)

Next we define hamiltonian vector fields and put a Lie algebra
structure on a subspace of
$C^{\infty}(M)$.
\begin{definition}\label{dirhamvf}
A function $f$ on a Dirac manifold $(M,L)$ is {\bf admissible}
 if there
exists a smooth vector field $X_f$ such that $X_f\oplus df$ is a
section of $L$. A vector field $X_f$ as above is called a
hamiltonian vector field of $f$. The set of admissible functions
forms a subspace $C^{\infty}_{adm}(M)$ of $C^\infty(M)$.
\end{definition}
 If $f$ is admissible then $df|_{L \cap TM}=0$. The converse
holds where the characteristic distribution $L \cap
TM$ has constant rank, but not in general.  In other words, $df$ can be
contained in $ \rho_{T^*M}(L)$ without being the image of a smooth
section of $L$; see Example \ref{notadm}.   Since any two
hamiltonian vector fields of an admissible function $f$ differ by
a characteristic vector field, which
annihilates any other
admissible function, we can make the following definition.

\begin{definition} The bracket on $C^{\infty}_{adm}(M)$ is given by
$\{f,g\}=X_g\cdot f$.
\end{definition}
This bracket differs by a sign from  the one in the original paper
of Courant \cite{Co}, but it allows us to recover the usual
conventions for presymplectic and Poisson manifolds, as shown
below.
 The main feature of this bracket is the following
(see \cite{Co}, Prop. 2.5.3):
\begin{proposition}\label{DiracLieAlg}
Let $(M,L)$ be a Dirac manifold. If $X_f$ and $X_g$ are any
hamiltonian vector fields for the admissible functions $f$ and
$g$, then $-[X_f,X_g]$ is a hamiltonian vector field for
$\{f,g\}$, which is therefore admissible as
  well.  The integrability of $L$ implies that
  the bracket satisfies the Jacobi identity, so
  $(C^{\infty}_{adm}(M),\{\cdot,\cdot\})$ is a Lie algebra.
\end{proposition}
We remark that the above can be partially extended to the space
$C^{\infty}_{bas}(M)$ of {\bf basic} functions, i.e. of functions
$\phi$ satisfying $d\phi|_{L \cap TM}=0$, which contains the
admissible functions. (This two spaces of functions coincide when
$L\cap TM$ is regular). Indeed, if $h$ is admissible and $\phi$ is
basic, then $\{\phi,h\}:=X_h\cdot \phi$ is well defined and basic,
since the flow of a hamiltonian vector field $X_h$ induces vector
bundle automorphisms of $TM\oplus T^*M$ that preserve $L\cap TM$
(see Section 2.4 in \cite{Co}). If $f$ is an admissible function,
then the Jacobiator of $f,h,$ and $\phi$ vanishes (adapt the proof
of
Prop. 2.5.3 in \cite{Co}).\\

We recall how manifolds endowed with 2-forms or bivectors fit into the
 framework of Dirac geometry. Let $\omega$ be a 2-form on $M$,
 $\tilde{\omega}:TM\rightarrow T^*M$ the bundle map $X \mapsto
 \omega(X,\cdot)$. Its graph $L=\{X\oplus\tilde{\omega}(X):X\in TM\}$
 is an almost Dirac structure; it is integrable iff $\omega$ is
 closed. If $\omega$ is symplectic, i.e. nondegenerate, then every
 function $f$ is admissible and has a unique hamiltonian vector field
 satisfying $\tilde{\omega}(X_f)=df$; the bracket is given by
 $\{f,g\}=\omega(X_f,X_g)$.

\begin{example} \label{notadm}
Let $\omega$ be the presymplectic form $x_1^2 dx_1 \wedge dx_2$ on
$M=\RR^2$, and let $L$ be its graph. The characteristic
distribution $L\cap TM$ has rank zero everywhere except along
$\{x_1=0\}$, where it has rank two, and it is clearly not
integrable (compare  the discussion following Definition
\ref{defdirstr}).  The differential of $f=x_1^2$
 takes all its values in the range of $\rho_{T^*M}$, but $f$
is  not admissible. This illustrates the remark
  following Definition
  \ref{dirhamvf}, i.e. it provides an example of a function which
  is basic but not admissible.
\end{example}

Let $\Lambda$ be a bivector field on $M$, $\tilde{\Lambda}:T^*M
\rightarrow TM$ the corresponding bundle map
$\xi\mapsto\Lambda(\cdot,\xi)$. (Note that the argument $\xi$ is
in the second position.)  Its graph
$L=\{\tilde{\Lambda}(\xi)\oplus \xi: \xi \in T^*M\}$ is an almost
Dirac structure which is integrable iff $\Lambda$ is a Poisson
bivector (i.e. the Schouten bracket $[\Lambda,\Lambda]_S$ is
zero). Every  function $f$ is admissible with a unique hamiltonian
vector field $X_f=\{\cdot,f\}$, and the bracket of functions is
$\{f,g\}=\Lambda(df,dg)$.

\section{Dirac-Jacobi manifolds} \label{secjacdir}

Dirac-Jacobi structures were introduced by Wade \cite{Wa} (under a different name)
and
include Jacobi (in particular, contact) and Dirac structures as
special cases.  Like Dirac structures, they are defined as maximal
isotropic subbundles of a certain vector bundle.

\begin{definition} A Dirac-Jacobi structure on a vector space $V$ is a
subspace $\bar{L}\subset (V\times \RR)\oplus (V^*\times \RR)$
which is maximal isotropic under the symmetric pairing
\begin{equation}\begin{split} \big \langle (X_1, f_1)\oplus(\xi_1, g_ 1)\;, \;(X_2,
f_2)\oplus(\xi_2, g_ 2) \big \rangle_+ = \frac{1}{2}(i_{X_2} \xi_1
+ i_{X_1} \xi_2 +g_1f_2+g_2f_1). \end{split}\end{equation}
\end{definition}
A Dirac-Jacobi structure on $V$ necessarily satisfies $\dim
\bar{L}=\dim V+1$. Furthermore, Equations \eqref{annih} hold for
Dirac-Jacobi structures too:
\begin{eqnarray}\rho_V(\bar{L})=(\bar{L} \cap V^*)^{\circ}
 \text{     and     }\rho_{V^*}(\bar{L})=(\bar{L} \cap
 V)^{\circ}.
\end{eqnarray}
As in the Dirac case, one has notions of pushforward and pullback
structures and as well as forward and backward maps. For example,
given a Dirac-Jacobi structure $\bar{L}$ on $V$ and a linear map
 $p: V \rightarrow Z$, one obtains
 a pushforward  Dirac-Jacobi structure on $Z$ by
$\{(pX,f)\oplus (\xi,g):(X,f)\oplus (p^*\xi,g)\in \bar{L}\}$.

On a manifold $M$, a maximal isotropic subbundle $\bar{L}\subset
\cE^1(M):= (TM\times \RR)\oplus (T^*M\times \RR)$ is called an
{\bf almost Dirac-Jacobi structure} on $M$.
\begin{definition}[\cite{Wa}, Def. 3.2]\label{defjdbr}
A {\bf Dirac-Jacobi structure}  on a manifold $M$ is an almost
Dirac-Jacobi structure $\bar{L}$ on $M$ whose space of sections is
closed under the extended Courant bracket on sections of
$\cE^1(M)$, which is defined by
\begin{equation}\label{extcoubra}
\begin{split}
[(X_1,f_1)\oplus(\xi_1,g_1)\;,\;(X_2,f_2)&\oplus(\xi_2,g_2)]=
\big([X_1,X_2],X_1\cdot f_2-X_2\cdot f_1\big)\\
&\oplus \big(\cL_{X_1}\xi_2-\cL_{X_2}\xi_1
+\frac{1}{2}d(i_{X_2}\xi_1
-i_{X_1}\xi_2)\\
&+
f_1\xi_2-f_2\xi_1+\frac{1}{2}(g_2df_1-g_1df_2-f_1dg_2+f_2dg_1),\\
&X_1\cdot g_2-X_2\cdot
g_1+\frac{1}{2}(i_{X_2}\xi_1-i_{X_1}\xi_2-f_2g_1+f_1g_2)\big).
\end{split}\end{equation}
\end{definition}

By a straightforward computation (see also Section 4 of
\cite{GrM}) this bracket can be derived from the
Courant bracket \eqref{coubra}, as follows.
 Denote by $U$ the
embedding  $\Gamma(\cE^1(M))\rightarrow \Gamma(T(M\times
\RR)\oplus T^*(M\times \RR))$ given by $$(X,f)\oplus(\xi,g)
\mapsto (X+f \frac{\partial }{\partial t})\oplus e^t(\xi+g dt),$$
where $t$ is the coordinate on the $\RR$ factor of the manifold
$M\times \RR$. Then $U$ is a bracket-preserving map from
$\Gamma(\cE^1(M))$ with the extended bracket \eqref{extcoubra}, to
$\Gamma(T(M\times \RR)\oplus T^*(M\times \RR))$ with the Courant
bracket \eqref{coubra} of the manifold $M\times \RR$.

Furthermore in Section 5 of \cite{IM} it is shown that any
Dirac-Jacobi manifold $(M,\bar{L})$ gives rise to a Dirac
structure on $M\times \RR$ given by
$$\tilde{\bar{L}}_{(x,t)}=\{(X+f\frac{\partial }{\partial
t})\;\oplus\;e^t(\xi+g dt): (X,f)\oplus(\xi,g)\in \bar{L}_x \},$$
where $t$ is the coordinate on $\RR$. This procedure extends the
well known symplectization of contact manifolds and Poissonization
of Jacobi manifolds, and may be called ``Diracization''.

If $\bar{L}$ is a Dirac-Jacobi structure,
$(\bar{L},\rho_{TM},[\cdot,\cdot])$ is a Lie algebroid (\cite{Wa},
Thm. 3.4), and each leaf of the induced foliation on $M$ has
 the structure of a precontact manifold (i.e. simply a 1-form)
or of a locally conformal presymplectic manifold (i.e. a 2-form
$\Omega$ and a closed 1-form $\omega$ satisfying
$d\Omega=\omega\wedge \Omega$). See Section \ref{leafDJ} for a
description of the induced foliation.
 As in the Dirac case, one can define
hamiltonian vector fields and endow a subset of $C^{\infty}(M)$
with a Lie algebra structure.
\begin{definition}[\cite{Wa}, Def. 5.1]
A function $f$ on a Dirac-Jacobi manifold $(M,\bar{L})$ is {\bf
admissible} if there exists a smooth vector field $X_f$ and a
smooth function $\varphi_f$ such that
$(X_f,\varphi_f)\oplus(df,f)$ is a section of $\bar{L}$. Pairs
$(X_f,\varphi_f)$ as above are unique up to smooth sections of
$\bar{L}\cap(TM \times \RR)$, and $X_f$ is called a hamiltonian
vector field of $f$. The set of admissible functions is denoted by
$C^{\infty}_{adm}(M)$.
\end{definition}

\begin{definition} The bracket on $C^{\infty}_{adm}(M)$ is given by
$\{f,g\}=X_g\cdot f+f\varphi_g$
\end{definition}

This bracket, which  differs by a sign from that
in \cite{Wa}, enjoys the
same properties stated in Proposition \ref{DiracLieAlg} for Dirac
manifolds (see \cite{Wa}, Prop. 5.2 and Lemma 5.3).

\begin{proposition}\label{Hamvf}
Let $(M,\bar{L})$ be a Dirac-Jacobi manifold. If $f$ and $g$ are
admissible functions, then
\begin{equation}
\begin{split}
&[(X_f,\varphi_f)\oplus(df,f)\;,\;(X_g,\varphi_g)\oplus(dg,g)]=\\
&([X_f,X_g],X_f\cdot \varphi_g-X_g\cdot
\varphi_f)\oplus(-d\{f,g\},-\{f,g\}),
\end{split}
\end{equation}
hence $\{f,g\}$ is again admissible. The integrability of
$\bar{L}$ implies that the admissible functions form a Lie
algebra.
\end{proposition}
We call a
function $\psi$ on $M$ {\bf basic} if  $X\cdot \psi+\psi f=0$ for
all elements $(X,f)\in \bar{L}\cap(TM\times\RR)$. This is
equivalent to requiring $(d\psi,\psi)\in \rho_{T^*M\times
\RR}(\bar{L})$ at each point of $M$, so the basic functions
contain the admissible ones. As in the case of Dirac structures,
we have the following properties:
\begin{lemma}\label{basicDJ}
 If $\psi$ is a basic and $h$ an admissible function, then the bracket
$\{\psi,h\}:=X_h\cdot \psi+\psi h$ is well-defined and again
basic.
\end{lemma}
\begin{proof}
It is clear that the bracket is well defined. To show that
$X_h\cdot \psi+\psi h$ is again basic we reduce the problem to the
Dirac case. Let
 $(X,f)\in
\bar{L}_x\cap(TM\times\RR)$ Fix a choice of $(X_h,\varphi_h)$ for
the admissible function $h$. The vector field
$X_h+\varphi_h\frac{\partial}{\partial t}$ on the Diracization
$(M\times \RR,\tilde{\bar{L}})$  (which is just a Hamiltonian
vector field of $e^th$) has a flow $\tilde{\phi}_{\epsilon}$,
which projects to the flow $\phi_{\epsilon}$ of $X_h$ under
$pr_1:M\times\RR\rightarrow M$. For each $\epsilon$ the flow
$\tilde{\phi}_{\epsilon}$ induces a vector bundle automorphism
$\Phi_{\epsilon}$ of $\cE^1(M)$, covering the diffeomorphism
$\phi_{\epsilon}$ of $M$, as follows:
\begin{eqnarray*}
(X,f)\oplus (\xi,g)\in \cE_x^1(M) \mapsto
(\tilde{\phi}_{\epsilon})_*(X\oplus f \frac{\partial}{\partial
t})_{(x,0)} \oplus
 (\tilde{\phi}_{\epsilon}^{-1})^*(\xi+g dt)_{(x,0)}\cdot
 e^{-pr_2(\tilde{\phi}_{\epsilon}(x,0))},
\end{eqnarray*}
where we identify $T_{\tilde{\phi}_{\epsilon}(x,0)}(M\times\RR)
\oplus T_{\tilde{\phi}_{\epsilon}(x,0)}^*(M\times\RR)$ with
$\cE^1_{\phi_{\epsilon}(x)}(M)$ to make sense of the second term.
Since the vector bundle maps induced by the flow
$\tilde{\phi}_{\epsilon}$ preserve the Dirac structure
$\tilde{\bar{L}}$ (see Section 2.4 in \cite{Co}), using the
definition of the Diracization $\tilde{\bar{L}}$ one sees that
$\Phi_{\epsilon}$ preserves $\bar{L}$, and therefore also
$\bar{L}\cap(TM\times \RR)$.
 Notice that we can pull back sections of $\cE^1(M)$ by
setting $(\Phi_{\epsilon}^*((X,f)\oplus(\xi,g)))_x:=
\Phi_{\epsilon}^{-1} ((X,f)\oplus(\xi,g))_{\phi_{\epsilon}(x)}.$ A
computation shows that
$$(0,0)\oplus (d(X_h\cdot \psi+\varphi_h\psi)\,,\,
X_h\cdot \psi+\varphi_h\psi)=
 \frac{\partial}{\partial
\epsilon}\Big|_0\Phi_{\epsilon}^*((0,0)\oplus(d\psi,\psi)),$$ so
that
\begin{equation*}\begin{split}
&\langle(0,0)\oplus (d(X_h\cdot
\psi+\varphi_h\psi),X_h\cdot \psi+\varphi_h\psi)\;,\;(X,f)\oplus(0,0)\rangle_+ =\\
&\frac{\partial}{\partial \epsilon}\Big|_0 \big[\langle
(0,0)\oplus(d\psi,\psi)_{\phi_{\epsilon}(x)}\;,\; \Phi_{\epsilon}
((X,f)\oplus(0,0))_x\rangle_+
e^{pr_2(\tilde{\phi}_{\epsilon}(x,0))}\big]=0,
\end{split}\end{equation*} as was to be shown.
\end{proof}
Furthermore, the Jacobiator of admissible functions $f,h$ and a
basic function $\psi$ is zero. One can indeed check that Wade's
proof of the Jacobi identity for admissible functions (\cite{Wa}
Prop. 5.2) applies in this case too. Alternatively, this follows
from the analogous statement for the Diracization $M\times \RR$,
since the map
\begin{eqnarray}\label{homJDiracDirac}  C^{\infty}_{adm}(M)
\rightarrow C^{\infty}_{adm}(M\times \RR)\;,\; g \mapsto e^tg
\end{eqnarray} is a well-defined Lie algebra
homomorphism\footnote{For the well-definedness notice that, if
$(X_g,\varphi_g)\oplus(dg,g)\in \Gamma(\bar{L})$, then
$(X_g+\varphi_g\frac{\partial }{\partial t})\oplus d(e^tg) \in
\Gamma(\tilde{\bar{L}})$. Notice that in particular
$X_g+\varphi_g\frac{\partial }{\partial t}$ is a hamiltonian
vector field for $e^tg$. Using this, the equation $e^t\{f,g\}_{M}=
\{e^tf,e^tg\}_{M\times \RR}$ follows at once from the definitions
of the respective brackets of functions.} and maps basic functions
to basic functions.\\

 Now we display some
examples of Dirac-Jacobi manifolds.

There is a one-to-one correspondence between Dirac structures on
$M$ and Dirac-Jacobi structures on $M$ contained in
$TM\oplus(T^*M\times \RR)$: to each Dirac structure $L$ one
associates the Dirac-Jacobi structure
$\{(X,0)\oplus(\xi,g):X\oplus \xi\in L, g\in \RR\}$ (\cite{Wa},
Remark 3.1).

A Jacobi structure on a manifold $M$ is given by a bivector field
$\Lambda$ and a vector field $E$ satisfying the Schouten bracket
conditions $[E,\Lambda]_S=0$ and $[\Lambda,\Lambda]_S=2E\wedge
\Lambda$.  When $E=0$, the Jacobi structure is a Poisson structure.
Any skew-symmetric vector bundle morphism $T^*M\times \RR
\rightarrow TM \times \RR$ is of the form $\left(
\begin{smallmatrix} \tilde{\Lambda} & -E \\ E & 0
\end{smallmatrix} \right)$ for a bivector field $\Lambda$ and a
vector field $E$, where as in Section \ref{secdir} we have
$\tilde{\Lambda}\xi=\Lambda(\cdot,\xi)$. $\text{Graph} \left(
\begin{smallmatrix} \tilde{\Lambda} & -E \\ E & 0
\end{smallmatrix} \right) \subset \cE^1(M)$ is a Dirac-Jacobi
structure iff $(M,\Lambda,E)$ is a Jacobi manifold (\cite{Wa},
Sect 4.1). In this case all functions are admissible, the unique
hamiltonian vector field of $f$ is\footnote{Again, this is opposite to
the usual sign convention.}
$X_f=\tilde{\Lambda} df- fE$,
 $\varphi_f = E\cdot f$ and the bracket is given by $\{f,g\}=\Lambda(df,dg)+fE\cdot g
-gE\cdot f$.

Similarly (see \cite{Wa}, Sect. 4.3), any skew-symmetric vector
bundle morphism $TM\times \RR \rightarrow T^*M \times \RR$ is of
the form $\left(
\begin{smallmatrix} \tilde{\Omega} & \sigma \\ -\sigma & 0
\end{smallmatrix} \right)$ for a 2-form $\Omega$ and a 1-form $\sigma$,
and  $\text{graph} \left(
\begin{smallmatrix} \tilde{\Omega} & \sigma \\ -\sigma & 0
\end{smallmatrix} \right) \subset \cE^1(M)$
is a Dirac-Jacobi structure iff $\Omega=d\sigma$.

Any contact form $\sigma$ on a manifold $M$ defines a Jacobi
structure $(\Lambda,E)$ (where $E$ is the Reeb vector field of
$\sigma$ and $\tilde{\Lambda}\tilde{d\sigma}|_{\ker
\sigma}=\text{Id}$; see for example \cite{IW}, Sect. 2.2), and
$\text{graph} \left(
\begin{smallmatrix} \tilde{d\sigma} & \sigma \\ -\sigma & 0
\end{smallmatrix} \right)$ is equal to $\text{graph} \left(
\begin{smallmatrix} \tilde{\Lambda} & -E \\ E & 0
\end{smallmatrix} \right)$. Further,
 by considering suitably defined graphs, one sees
that locally conformal presymplectic structures and homogeneous
Poisson manifolds (given by a Poisson bivector $\Lambda$ and a
vector field $Z$ satisfying $\cL_Z\Lambda=-\Lambda$) are examples
of Dirac-Jacobi structures (\cite{Wa}, Sect. 4).

\section{The prequantization spaces}\label{spaces}
In this section  we determine the prequantization condition for a
Dirac manifold $(P,L)$, and we describe its ``prequantization space''
(i.e. the geometric object that allows us to find a representation
of $C^{\infty}_{adm}(P)$).

 We recall the prequantization of
 a Poisson manifold $(P,\Lambda)$ by a $U(1)$-bundle as described
in \cite{CMdL}. The bundle map $\tilde{\Lambda}:T^*P\rightarrow
TP$ extends to a cochain map from forms to multivector fields,
which descends to a map from de Rham cohomology
$H^{\bullet}_{dR}(P,\RR)$ to Poisson cohomology
$H^{\bullet}_{\Lambda}(P)$ (the latter having the set of
$p$-vector fields as $p$-cochains). The prequantization condition,
first formulated in this form in \cite{Va2},
 is
that $[\Lambda]\in H^{2}_{\Lambda}(P)$ be the image under
$\tilde{\Lambda}$ of an integral de Rham class, or equivalently
that
\begin{eqnarray}\label{preqpois}
\tilde{\Lambda} \Omega=\Lambda+{\cL}_A\Lambda \end{eqnarray} for
some integral closed 2-form $\Omega$ and vector field $A$ on $P$.
Assuming this prequantization condition to be satisfied, let
$\pi:Q\rightarrow P$ be a $U(1)$-bundle with first Chern class
$[\Omega]$, $\sigma$ a connection on $Q$ with curvature $\Omega$
(i.e. $d\sigma=\pi^*\Omega$), and $E$ the generator of the
$U(1)$-action (so that $\sigma(E)=1$ and  $\pi_*E=0$). Then (see
Thm. 3.1 in \cite{CMdL})
\begin{eqnarray}\label{jac} (\Lambda^H+E\wedge A^H,E)
\end{eqnarray}
is a Jacobi structure on $Q$ which pushes down to $(\Lambda,0)$ on $P$
via $\pi_*$.
(The superscript
 $^H$ denotes horizontal lift, with respect to the connection
$\sigma$, of multivector fields on $P$.)  We say that $\pi$ is a
 Jacobi map.

It follows from the Jacobi map property of $\pi$ that assigning to
a function $f$ on $P$ the  hamiltonian vector field of
$-\pi^*f$, which is
 $-\widetilde{(\Lambda^H+E\wedge A^H)}
(\pi^*df)+(\pi^*f)E$, defines a Lie algebra
homomorphism from $C^\infty(P)$ to the operators on $C^\infty(Q)$.\\

 Now we carry out an analogous construction on a Dirac manifold
$(P,L)$.  Recall
that $L$ is a Lie algebroid with the restricted Courant bracket and
 anchor $\rho_{TP}: L
\rightarrow TP$ (which is just the projection onto the tangent
component).  This anchor gives a Lie algebra homomorphism from
$\Gamma(L)$  to $\Gamma(TP)$ with the Lie bracket of vector
fields.  The pullback by the anchor therefore induces a map
$\rho_{TP}^*: \Omega^{\bullet}_{dR}(P,\RR)\rightarrow
\Omega^{\bullet}_{L}(P)$, descending to a map from de Rham
cohomology to the Lie algebroid cohomology $H^2_L(P)$.  (We recall
from \cite{CW} that $\Omega^{\bullet}_{L}(P)$ is the graded
differential algebra of sections of the exterior algebra of
$L^*$.)
 There is a distinguished class in $H^{2}_{L}(P)$:
on $TP\oplus T^*P$, in addition to the natural symmetric pairing
\eqref{symm}, there is also an anti-symmetric one given by
\begin{eqnarray}
\langle X_1\oplus \xi_1,X_2\oplus \xi_2\rangle_- = \frac{1}{2}(
i_{X_2}\xi_1- i_{X_1}\xi_2).
\end{eqnarray}
Its restriction $\Upsilon$ to $L$ satisfies $d_L \Upsilon=0$. Our
prequantization condition is
\begin{eqnarray}\label{cond0} [\Upsilon] \in
\rho_{TP}^*(i_*(H^2(P,\ZZ)))
\end{eqnarray}
or equivalently
\begin{eqnarray} \label{cond1}\rho_{TP}^*\Omega=\Upsilon+d_L \beta,
\end{eqnarray}
where $\Omega$ is a closed integral 2-form and $\beta$ a 1-cochain
for the Lie algebroid $L$, i.e. a section of $L^*$.

\begin{remark}\label{remcond} If $L$ is the graph of a
presymplectic form $\omega$ then $\Upsilon = \rho_{TP}^*(\omega) $.
If $L$ is $\text{graph} (\tilde{\Lambda})$ for a
Poisson bivector $\Lambda$ and $\Omega$ is a 2-form, then
$\rho_{TP}^*[\Omega]=[\Upsilon]$ if and only if
$\tilde{\Lambda}[\Omega]=[\Lambda]$.\footnote{This
is consistent with the fact that, if $\omega$
is symplectic, then \text{graph}$(\tilde{\omega})=
\text{graph}(\tilde{\Lambda})$, where the bivector $\Lambda$  is
defined so that the vector bundle maps $\tilde{\omega}$ and
$\tilde{\Lambda}$ are inverses of each other (so if
$\omega=dx\wedge dy$ on $\RR^2$, then
$\Lambda=\frac{\partial}{\partial x} \wedge
\frac{\partial}{\partial y}$).} $[\Omega]$ to $[\Upsilon]$. This
shows that \eqref{cond0} generalizes the
prequantization conditions for presymplectic and Poisson
structures mentioned in the introduction and in formula
\eqref{preqpois}.
\end{remark}

\begin{remark}\label{twisted}
The prequantization condition
above can not even be formulated for twisted Dirac structures.
We recall the definition of these structures \cite{SeW}.
If
 $\phi$ is a closed 3-form on a manifold $P$,
adding the term $\phi(X_1,X_2,\cdot)$ to  the Courant bracket
(i.e. to the right hand side of Equation \eqref{coubra})
determines a new bracket $[\cdot,\cdot]^{\phi}$ so that $TP\oplus
T^*P$, together with this bracket, the original anchor $\rho_{TP}$
and the symmetric pairing $\langle\cdot,\cdot\rangle_+$, form a
Courant algebroid.  A $\phi$-twisted Dirac structure $L$ is then a
maximal isotropic subbundle which is closed under
$[\cdot,\cdot]^{\phi}$; it is automatically a Lie algebroid
(whose Lie algebroid differential we denote by $d_L^{\phi}$). The
orbits of the Lie algebroid carry 2-forms $\Omega_L$
given as in the remark following Definition \ref{lindirstr},
satisfying $d\Omega_L=j^*\phi$ where $j$ is the inclusion of a
leaf in $P$ and $d$ is the de Rham differential on the leaves.
Since
$$d_L^{\phi}\Upsilon=d_L^{\phi}\rho_{TP}^*\Omega_L= \rho_{TP}^*d
\Omega_L=\rho_{TP}^*j^*\phi,$$ we conclude that $\Upsilon$ is
usually not $d_L^{\phi}$-closed, so we cannot expect $\Omega$ to
be closed in  \eqref{cond1}, and hence we cannot require that it be
integral.
The correct notion of prequantization should probably involve a gerbe.
\end{remark}

Now, assuming the prequantization condition \eqref{cond0} and proceeding
as in the Poisson case,  let  $\pi:Q\rightarrow P$ be a
$U(1)$-bundle with connection form $\sigma$
having curvature $\Omega$; denote by $E$ the
infinitesimal generator of the $U(1)$-action.

\begin{theorem} \label{thmpreq}
The subbundle $\bar{L}$ of $\cE^1(Q)$ given by the direct
sum of
$$\{(X^H+\langle X\oplus \xi, \beta \rangle
E,0)\oplus(\pi^*\xi,0):
 X\oplus \xi\in L\}$$ and the line bundles generated by
 $(-E,0)\oplus(0,1)$ and
 $(-A^H,1)\oplus(\sigma-\pi^*\alpha, 0)$
is a Dirac-Jacobi structure on $Q$. Here, $A\oplus \alpha$ is an
isotropic section of $TP\oplus T^*P$ satisfying $\beta=2\langle
A\oplus \alpha, \,\cdot\, \rangle_+|_L$. Such a section always
exists, and the subbundle above is independent of the choice of
$A\oplus \alpha$.
\end{theorem}

\begin{proof}
Let $C$ be a maximal isotropic (with respect to $\langle \cdot,
\cdot \rangle_+$) complement of $L$ in $TP\oplus T^*P$. Such a
complement always exists, since the space of complements at each point
is contractible (an affine space modeled on a space of skew-symmetric
forms).
Now extend $\beta$ to a functional
$\tilde{\beta}$ on $TP\oplus T^*P$ by setting
$\tilde{\beta}|_C=0$. There exists a unique section $A\oplus
\alpha$ of $TP\oplus T^*P$ satisfying $\beta=2\langle A\oplus
\alpha, \,\cdot \,\rangle_+$ since the symmetric pairing is
non-degenerate. Since $\langle A\oplus \alpha, \,\cdot
\,\rangle_+|_C=0$ and $C$ is maximal isotropic we conclude that
$A\oplus \alpha$ belongs to $C$ and is hence isotropic itself.
This shows the
  existence of $A\oplus \alpha$ as above.

Now clearly $A\oplus \alpha+Y\oplus \eta$ satisfies the property
stated in the theorem iff $Y\oplus \eta\subset L$, and in this
case it is isotropic (i.e. $\langle A+Y,\alpha+\eta \rangle=0$)
iff $Y\oplus \eta\subset \ker \beta$. So a section $A\oplus
\alpha$ as in the theorem is unique up to sections $Y\oplus \eta$
of $\ker \beta$. By inspection one sees that replacing $A\oplus
\alpha$ by $A\oplus\alpha+Y\oplus\eta$ in the formula for
$\bar{L}$ defines the same subbundle.

That $\bar{L}$ is isotropic with respect to the symmetric pairing
on $\cE^1(Q)$ follows from the fact that  $L$ is isotropic,
together with the properties of $A\oplus \alpha$. $\bar{L}$ is
clearly a subbundle of dimension $\dim P+2$, so it is an almost
Dirac-Jacobi structure.

To show that $\bar{L}$ is integrable, we use  the  fact that
$\bar{L}$ is integrable if and only if $\langle [e_1,e_2],e_3
\rangle_+=0$ for all sections $e_i$ of $\bar{L}$ and that $
\langle [\cdot,\cdot],\cdot \rangle_+$ is a totally skew-symmetric
tensor if restricted to sections of $\bar{L}$, i.e. an element of
$\Gamma(\wedge^3 \bar{L}^*)$ (\cite{IM}, Prop. 2.2). Each section
of $\bL$ can be written as a $C^{\infty}(Q)$-linear combination of
the following three types of sections of $\bar{L}$:
$a=:(X^H+\langle X\oplus \xi, \beta \rangle
E,0)\oplus(\pi^*\xi,0)$ where $X\oplus \xi\in \Gamma(L)$,
$b:=(-E,0)\oplus(0,1)$ and $c:=(-A^H,1)\oplus(\sigma-\pi^*\alpha,
0)$.
We will use subscripts to label more than one section of a given type.
It is immediate that brackets of the form $[a,b]$,$[b_1,b_2]$,
and $[c_1,c_2]$ all vanish, and a computation shows that $ \langle
[a_1,a_2],a_3 \rangle_+=0$ since $L\subset TP\oplus T^*P$ is a
Dirac structure. Finally $ \langle [a_1,a_2],c \rangle_+=0$ using
$d\sigma=\pi^*\Omega$ and the prequantization condition
\eqref{cond1}, which when applied to sections $X_1\oplus \xi_1$
and $X_2\oplus \xi_2$ of $L$ reads
\begin{equation*}
 \textbf{} \Omega(X_1,X_2)=\langle \xi_1, X_2
\rangle +X_1 \langle \beta,X_2\oplus \xi_2\rangle -X_2 \langle
\beta, X_1\oplus \xi_1\rangle -
 \langle \beta\;,\;
\big[ X_1\oplus \xi_1,X_2\oplus \xi_2\big] \rangle.
\end{equation*}
By skew-symmetry, the vanishing of these expressions is enough to
prove the integrability of $\bL$.
\end{proof}

\begin{remark}
When $(P,L)$ is a Poisson manifold, $\bar{L}$ is exactly the graph
of the Jacobi structure \eqref{jac}, i.e. it generalizes the
construction of \cite{CMdL}. If $(P,L)$ is given by a
presymplectic form $\Omega$, then $\bar{L}$ is the graph of
$(d\sigma,\sigma)$.
\end{remark}

\begin{remark}  The construction of
Theorem \ref{thmpreq} also works for complex Dirac structures
(i.e., integrable maximal isotropic complex subbundles of the
complexified bundle $T_\CC M \oplus T_\CC^* M$). It can be adapted
to the setting of generalized complex structures \cite{Gua}
(complex Dirac structures which are transverse to their complex
conjugate) and generalized contact structures \cite{IW} (complex
Dirac-Jacobi structures which are transverse to their complex
conjugate) as follows. If $(P,L)$ is a generalized complex
manifold, assume all of the previous notation and the following
prequantization condition:
\begin{eqnarray} \rho_{TP}^*\Omega=i\Upsilon+d_L \beta,
\end{eqnarray}
where $\Omega$ is (the complexification of) a closed integer
2-form and $\beta$ a 1-cochain for the Lie algebroid $L$. Then the
direct sum of
$$\{(X^H+\langle X\oplus \xi, \beta \rangle
E,0)\oplus(\pi^*\xi,0):
 X\oplus \xi\in L\}$$
and the complex line bundles generated by $(-iE,0)\oplus(0,1)$
 and
$(-A^H,i)\oplus(\sigma-\pi^*\alpha, 0)$
is a generalized contact structure on $Q$, where $A\oplus
\alpha$ is the unique section of the conjugate of $L$ satisfying
$\beta=2\langle A\oplus \alpha, \,\cdot\, \rangle_+|_L$.
\end{remark}

\subsection{Leaves of the Dirac-Jacobi structure}\label{leafDJ}
 Given any Dirac-Jacobi manifold
$(M,\bar{L})$, each leaf of the foliation integrating the
distribution $\rho_{TM}(\bar{L})$ carries one of two kinds of
geometric structures \cite{IM}, as we describe now.
$\rho_1:\bar{L}\rightarrow \RR, (X,f)\oplus(\xi,g) \mapsto f$
determines an algebroid 1-cocycle, and a leaf $\bar{F}$ of the
foliation will be of one kind or the other depending on whether
$\ker \rho_1$ is contained in the kernel of the anchor $\rho_{TM}$
or not.   (This property is satisfied either  at all points of
$\bar{F}$ or at none). As with Dirac structures, the Dirac-Jacobi
structure $\bar{L}$ determines a field of skew-symmetric bilinear
forms $\Psi_{\bar{F}}$ on the image of $\rho_{TM}\times \rho_1$.

 If $\ker \rho_1
\not\subset \ker \rho_{TM}$ on $\bar{F}$ then $\rho_{TM}\times
\rho_1$ is surjective, hence $\Psi_{\bar{F}}$ determines a 2-form
and a 1-form on $\bar{F}$. The former is the differential of the
latter, so the leaf $\bar{F}$ is simply endowed with a 1-form ,
i.e. it is a precontact leaf. If $\ker \rho_1\subset \ker
\rho_{TM}$ on $\bar{F}$ then the image of $\rho_{TQ}\times \rho_1$
projects isomorphically onto $T\bar{F}$, which therefore carries
a 2-form $\Omega_{\bar{F}}$. It turns out that
$\omega_{\bar{F}}(Y):=-\rho_1(e)$, for any $e\in \bar{L}$ with
$\rho_{TM}(e)=Y$, is a well-defined 1-form on $\bar{F}$, and that
$(\bar{F},\Omega_{\bar{F}},\omega_{\bar{F}})$ is a locally
conformal presymplectic manifold, i.e. $\omega_{\bar{F}}$ is
closed and
$d\Omega_{\bar{F}}=\Omega_{\bar{F}}\wedge\omega_{\bar{F}}$.

On our prequantization $(Q,\bar{L})$
the leaf $\bar{F}$ through $q\in Q$ will carry one or the other
geometric structure depending on whether $A$ is tangent to $F$,
where $F$ denotes the presymplectic leaf of $(P,L)$ passing
through $\pi(q)$. Indeed one can check that at $q$ we have $\ker
\rho_1 \not\subset \ker \rho_{TQ} \Leftrightarrow A\in
T_{\pi(q)}F$. When $\ker \rho_1 \not\subset \ker \rho_{TQ}$ on a
leaf $\bar{F}$ we hence deduce that $\bar{F}$, which is equal to
$\pi^{-1}(F)$, is a precontact manifold, and a computation shows
that the 1-form is given by the restriction of
$$\sigma+\pi^*(\xi_A-\alpha)$$
 where $\xi_A$
is any covector satisfying $A\oplus \xi_A \in L$.

A leaf $\bar{F}$ on which
 $\ker \rho_1 \subset \ker \rho_{TQ}$ is locally
 conformal presymplectic, and its image under $\pi$ is
an integral submanifold of the integrable distribution
$\rho_{TP}(L)\oplus\RR A$ (hence a one parameter family
 of presymplectic leaves).
  A computation shows that
the locally conformal presymplectic structure is given by

$$(\omega_{\bar{F}},
\Omega_{\bar{F}})=\left(\pi^*\tilde{\gamma}\;,\;
 (\sigma-\pi^*\alpha)\wedge\pi^*\tilde{\gamma}
 +\pi^*\tilde{\Omega}_L\right).$$
 Here $\tilde{\gamma}$ is the 1-form on $\pi(\bar{F})$
with kernel $\rho_{TP}(L)$ and evaluating to one on
 $A$, while $\tilde{\Omega}_L$ is the two form on $\pi(\bar{F})$ which
coincides with $\Omega_L$ (the presymplectic form on the leaves of
$(P,L)$) on $\rho_{TP}(L)$ and annihilates $A$.

\subsection{Dependence of the Dirac-Jacobi structure on
choices}\label{choices}

Let $(P,L)$ be a prequantizable Dirac manifold , i.e. one for
which there exist  a closed integral 2-form $\Omega$  and a
section of $\beta$ of $L^*$ such that
\begin{eqnarray}\label{againcond1}\rho_{TP}^*\Omega=\Upsilon+d_L \beta.\end{eqnarray}
 The Dirac-Jacobi manifold $(Q,\bar{L})$ as
defined in Theorem \ref{thmpreq} depends on three data: the choice
(up to isomorphism) of the $U(1)$-bundle $Q$, the choice of
connection $\sigma$ on $Q$
whose curvature has cohomology class $i_*c_1(Q)$, and the choice
of $\beta$, subject to the condition that Equation
\eqref{againcond1} be satisfied. We will explain here how
the Dirac-Jacobi structure $\bar{L}(Q,\sigma,\beta)$
 depends on these choices.

First, notice that the value of $\Omega$ outside of $\rho_{TP}(L)$
does not play a role in \eqref{againcond1}. In fact, different
choices of  $\sigma$ agreeing over $\rho_{TP}(L)$  give rise to the same
Dirac-Jacobi structure. This is consistent with the following
lemma, which is the result of a straightforward computation:
\begin{lemma} \label{lemmathreeNEW}For any 1-form $\gamma$ on $P$
the Dirac-Jacobi structures $\bar{L}(Q,\sigma,\beta)$ and
$\bar{L}(Q,\sigma+\pi^*\gamma,\beta+\rho_{TP}^*\gamma)$ are equal.
\end{lemma}

Two Dirac-Jacobi structures on a given $U(1)$-bundle $Q$ over $P$
give isomorphic quantizations if they are related by an element of
the gauge group $C^\infty(P,U(1))$ acting on $Q$.  Noting that the
Lie algebroid differential $d_L$ descends to a map
$C^\infty(P,U(1))\to \Omega^1_L(P)$ we denote by $H^1_L(P,U(1))$
the quotient of the closed elements of $\Omega^1_L(P)$ by the
space $d_L(C^\infty(P,U(1))$ of $U(1)$-exact forms.

Now  we  show:
\begin{proposition}\label{propchoice}
The set of isomorphism classes of Dirac-Jacobi manifolds prequantizing
$(P,L)$ maps surjectively to the space $(\rho_{TP}^* \circ
i_*)^{-1}[\Upsilon]$ of topological types of compatible $U(1)$-bundles; the prequantizations of a given topological type
are a principal homogeneous space for  $ H^1_L(P,U(1)).$
\end{proposition}
\begin{proof}
Make a choice of prequantizing triple $(Q, \sigma, \beta)$.  With
$Q$ and $\sigma$ fixed, we are allowed to change $\beta$ by a
$d_L$-closed section of $L^*$.  If we fix only $Q$, we are allowed
to change $\sigma$ in such a way that the resulting curvature
represents the cohomology class $i_*c_1(Q)$, so we can change
$\sigma$ by $\pi^*\gamma$ where $\gamma$ is a 1-form on $P$. Now
$\bar{L}(Q,\sigma+\pi^*\gamma,\tilde{\beta})=
\bar{L}(Q,\sigma,\tilde{\beta}-\rho_{TP}^*\gamma)$ by Lemma
\ref{lemmathreeNEW}, so we obtain one of the Dirac-Jacobi
structures already obtained above.  Now, if we replace $\beta$ by
 $\beta + d_L\phi$ for $\phi \in C^\infty(P,U(1))$,
we obtain an isomorphic Dirac-Jacobi structure: in fact
$\bar{L}(Q,\sigma,\beta)$ is equal to $\bar{L}(Q,\sigma+\pi^*d
\phi,\beta+d_L \phi)$ by Lemma \ref{lemmathreeNEW}, which is
isomorphic to $\bar{L}(Q,\sigma,\beta+d_L \phi)$ because the gauge
transformation given by $\phi$ takes the connection $\sigma$ to
$\sigma+\pi^*d \phi$.
 So we see that the difference between two
prequantizing Dirac-Jacobi structures on the fixed $U(1)$-bundle
$Q$ corresponds to an element of  $H^1_L(P,U(1))$.

\end{proof}

In Dirac geometry, a B-field transformation (see for example
 \cite{SeW}) is an automorphism of the Courant algebroid $TM\oplus
 T^{*}M$ arising from a closed 2-form $B$ and taking each Dirac
 structure into another one with an isomorphic Lie algebroid.  There
 is a similar construction for Dirac-Jacobi structures.  Given any
 1-form $\gamma$ on any manifold $M$, the vector bundle endomorphism
 of $\cE^1(M)= (TM\times \RR)\oplus (T^*M\times \RR)$ that acts on
 $(X,f)\oplus(\xi,g)$ by adding $(0,0)\oplus \left(
 \begin{smallmatrix} \tilde{d\gamma} & \gamma \\ -\gamma & 0
\end{smallmatrix} \right)(X,f)$ preserves the extended Courant
bracket and the symmetric pairing.  Thus, it maps each Dirac-Jacobi structure to another one.  We call this operation an extended B-field transformation.

\begin{lemma}\label{lemmatwoNEW}
Let $\gamma$ be a closed 1-form on $P$.
Then $\bar{L}(Q,\sigma+\pi^*\gamma,\beta)$ is obtained from
$\bar{L}(Q,\sigma,\beta)$ by the extended B-field transformation associated to $\gamma$.
\end{lemma}

In the statements that follow, until the end of this subsection,
we assume that the distribution $\rho_{TP}(L)$ has constant rank,
and we denote by $\cF$ the regular distribution integrating it.
\begin{corollary}
Assume that $\rho_{TP}(L)$ has
 constant rank. Then the isomorphism classes of
prequantizing Dirac-Jacobi structures on the fixed $U(1)$-bundle
$Q$, up to extended B-field transformations, form a principal
homogeneous space for
$$H^1_L(P,U(1))/H^1_{\rho_{TP}(L)}(P,U(1)),$$
 where $H^{\bullet}_{\rho_{TP}(L)}(P)$ denotes the foliated
  (i.e. tangential de Rham) cohomology
 of  $\rho_{TP}(L)$.
\end{corollary}
\begin{proof}
We saw  in the proof of Prop. \ref{propchoice} that, if $(P,L)$ is
prequantizable, the prequantizing Dirac-Jacobi structures on a
fixed $U(1)$-bundle $Q$ are given by
$\bar{L}(Q,\sigma,\beta+\beta')$ where $Q,\sigma,\beta$ are fixed
and $\beta'$ ranges over all $d_L$-closed sections of
$L^*$. Consider
$\rho_{TP}^*\gamma$ for a closed 1-form $\gamma$. Then
$\bar{L}(Q,\sigma,\beta+\rho_{TP}^*\gamma)=L(Q,\sigma-\pi^*\gamma,\beta)$
by Lemma \ref{lemmathreeNEW}, and this is  related to
$\bar{L}(Q,\sigma, \beta)$ by an extended B-field transformation
because of Lemma  \ref{lemmatwoNEW}. To finish the argument,
divide by the $U(1)$-exact forms.
\end{proof}
 We will
now give a characterization of the $\beta$'s appearing in a
prequantization triple.

\begin{lemma}\label{regular}
Let $(P,L)$ be a Dirac manifold for which  $\rho_{TP}(L)$ is a
regular foliation.  Given a section $\beta'$ of $L^*$, write
$\beta'=\langle A'\oplus\alpha',\cdot \rangle|_L$.  Then
 $d_L\beta'=\rho_{TP}^*\Omega'$
for some 2-form along  $\cF$ iff the vector field $A'$ preserves
the foliation $\cF$. In this case,
$\Omega'=d\alpha'-\cL_{A'}\Omega_L$ where $\Omega_L$ is the
presymplectic form on the leaves of  $\cF$ induced by $L$.
\end{lemma}
\begin{proof}
For all sections $X_i\oplus \xi_i$ of $L$ we have
\begin{eqnarray*}\label{Omega'}
d_L\beta'(X_1\oplus \xi_1,X_2\oplus
\xi_2)=d\alpha'(X_1,X_2)+(\cL_{A'}\xi_2)X_1-(\cL_{A'}\xi_1)X_2+A'\cdot\langle
\xi_1,X_2 \rangle.
\end{eqnarray*}
Clearly  $d_L\beta'$ is of  the form $\rho_{TP}^*\Omega'$ iff
$L\cap T^*P\subset \ker d_L\beta'$ (and in this case $\Omega'$ is
clearly unique). Using the constant rank assumption to extend
appropriately elements of $L\cap T^*P$ to some neighborhood in
$P$, one sees that this is equivalent to $(\cL_{A'}\xi)X=0$ for
all sections $\xi$ of $L\cap T^*P=(\rho_{TP}(L))^{\circ}$  and
vectors $X$ in $\rho_{TP}(L)$, i.e. to $A'$ preserving the
foliation.

The formula for $\Omega'$ follows from a computation manipulating
the above expression for $d_L\beta'$ by means of the Leibniz rule
for Lie derivatives.
\end{proof}

We saw  in the proof of Prop. \ref{propchoice} that, if $(P,L)$ is
prequantizable, the prequantizing Dirac-Jacobi structures on a
fixed $U(1)$-bundle $Q$ are given by
$\bar{L}(Q,\sigma,\beta+\beta')$ where $Q,\sigma,\beta$ are fixed
and $\beta'$ ranges over all $d_L$-closed sections of $L^*$.
 Since
$\Upsilon=\rho_{TP}^*\Omega_L$, it follows from \eqref{againcond1}
 that $d_L\beta$ is the pullback by $\rho_{TP}$ of some
2-form along $\cF$. So, by the above lemma, $\beta=\langle
A\oplus\alpha,\cdot \rangle|_L$ for some vector field $A$
preserving the regular foliation $\cF$. Also, $\beta'=\langle
A'\oplus\alpha',\cdot \rangle|_L$ where $A'$ is a vector field
preserving $\cF$ and $d\alpha'-\cL_{A'}\Omega_L=0$, and conversely
every $d_L$-closed $\beta'$  arises this way (but choices of
$A'\oplus \alpha'$ differing by sections of $L$ will give rise to
the same $\beta'$).

\begin{example}
Let $F$ be an integrable distribution on a manifold $P$ (tangent to
 a regular foliation $\cF$), and $L=F\oplus
F^{\circ}$ the corresponding Dirac structure. By Lemma
\ref{regular}
 (or by
a direct computation) one sees that the $d_L$-closed sections
$\beta$ of $L^*$ are sums of sections of $TP/F$ preserving the
foliation and closed 1-forms along $F$. By Prop. \ref{propchoice},
the set of isomorphism classes of prequantizing Dirac-Jacobi
 structures maps surjectively  to the set
$\ker(\rho_{TP}^* \circ i_*)$ of topological types; the inverse image of a given type is a principal homogeneous space for  $$\{ \text{Sections of $TP/F$ preserving the foliation}\}\times H^1_{F}(P,U(1)),$$
where the Lie algebroid cohomology
 $H^{\bullet}_F(P)$ is the tangential de Rham cohomology of $\cF$
(and $\ker(\rho_{TP}^* \circ i_*)$ denotes the kernel in degree
two).
\end{example}


\section{The prequantization representation}\label{representation}

In this section, assuming the prequantization condition
\eqref{cond1} for the Dirac manifold $(P,L)$ and denoting by
$(Q,\bar{L})$ its prequantization as in Theorem \ref{thmpreq}, we
construct a representation of the Lie algebra
$C^{\infty}_{adm}(P)$. We will do so by first mapping this space
of functions to a set of ``equivalence classes of vector fields''
on $Q$ and then by letting these act on
$C^{\infty}_{bas}(Q,\CC)_{P-loc}$, a sheaf over $P$. Here
$C^{\infty}_{bas}(Q,\CC)$ denotes the complex basic\footnote{We
use basic instead of admissible functions in order to obtain the
same representation as in Section \ref{line}.}
 functions on
$(Q,\bar{L})$, as defined in Section \ref{secjacdir}, which in the
case at hand are exactly the functions whose differentials
annihilate $\bar{L}\cap TQ$.
The subscript ``$_{P-loc}$'' indicates that we consider functions
which are defined on subsets $\pi^{-1}(U)$ of $Q$, where $U$
ranges over the open subsets\footnote{We use the space of
\emph{P-local} instead of \emph{global} basic functions because
the latter could be too small for certain injectivity statements.
See Proposition \ref{inj} below and the remarks following it, as
well as Section \ref{remarks}.} of $P$. We will decompose this
representation and make some comments on the faithfulness of the
resulting
subrepresentations.\\

Let  $\tilde{L}=\{(X,0)\oplus(\xi,g): X\oplus \xi\in L, g\in
\RR\}$ be the Dirac-Jacobi structure associated to the Dirac
structure $L$ on $P$. It is immediate that $\tilde{L}$ is the
push-forward of $\bar{L}$ via $\pi:Q\rightarrow P$, i.e.
$\tilde{L}=\{(\pi_*Y,f)\oplus(\xi,g):(Y,f)\oplus(\pi^*\xi,g)\in
\bar{L}\}.$
  From this it follows that if functions $f,g$ on $P$ are admissible
then their pullbacks $\pi^*f,\pi^*g$ are also
admissible\footnote{To show the smoothness of the hamiltonian
vector fields of $\pi^*f$ and $\pi^*g$, we actually have to use
the particular form of $\bar{L}$.}
 and
 \begin{eqnarray}
\label{pihom}\{\pi^*f,\pi^*g\}=\pi^*\{f,g\}.\end{eqnarray}

\begin{proposition}\label{therepr}
The map
\begin{equation}\begin{split}\label{homom}
(C^{\infty}_{adm}(P),\{\cdot,\cdot\}) &\rightarrow
\text{Der}(C^{\infty}_{bas}(Q,\CC)_{P-loc})\\
 g\;\;\;\;\;\;\;\;\; &\mapsto \{\pi^*g,\cdot\}
\end{split}\end{equation}
determines a representation of $C^{\infty}_{bas}(Q,\CC)_{P-loc}$.
\end{proposition}
\begin{proof}  Recall that the
expression $\{\pi^*g,\phi\}$ for $\phi\in
C^{\infty}_{bas}(Q,\CC)_{P-loc}$ was defined in Section
\ref{secjacdir} as $-X_{\pi^*g}(\phi)-\phi\cdot
0=-X_{\pi^*g}(\phi)$, for any choice $X_{\pi^*g}$
 of hamiltonian vector field for $\pi^*g$.
The proposition follows from the versions of the following
statements for basic functions (see Lemma \ref{basicDJ} and the
remark following it). First: the map \eqref{homom} is well-defined
since the set of admissible functions on the Dirac-Jacobi manifold
$Q$ is closed under the bracket $\{\cdot,\cdot\}$. Second:
 it is a Lie algebra homomorphism because of Equation \eqref{pihom} and
because the bracket of admissible functions on $Q$ satisfies the
Jacobi identity.
Alternatively, for the second statement we can make use of the
relation $[-X_{\pi^*f},-X_{\pi^*g}]=-X_{\{\pi^*f,\pi^*g\}}$ (see
Proposition \ref{Hamvf}).
\end{proof}

Since the
Dirac-Jacobi structure on $Q$ is invariant under the $U(1)$
action, the infinitesimal generator $E$ is a derivation of
the bracket. We can decompose $C^{\infty}_{bas}(Q,\CC)_{P-loc}$
into the eigenspaces $H^n_{bas}$ of $E$
  corresponding to the eigenvalues $2\pi i n$, where $n$ must be an
integer, and similarly for $H^n_{adm}$.  The derivation property
implies that
 $\{H^n_{adm},H^{n'}_{bas}\}\subseteq
H^{n+n'}_{bas}$. The Lie algebra of admissible functions on $P$
may be identified with the  real-valued global functions in
$H^0_{adm}$, which acts on each $H^n_{bas}$ by the bracket, i.e.
by the representation \eqref{homom}. In particular, the action on
 $H^{-1}_{bas}$
 is the
usual prequantization action.
 The classical limit is
obtained by letting $n\to -\infty$. Clearly all of the above
applies if we restrict the representation \eqref{homom} to
$C^{\infty}_{adm}(Q,\CC)_{P-loc}$, i.e. if we replace
``$H^n_{bas}$'' by ``$H^n_{adm}$'' above.
\\

Now we will comment on the faithfulness of the above
representations.  The map that assigns to an admissible function
$g$ on $P$ the equivalence class of hamiltonian vector fields of
$-\pi^*g$ depends on the choices of $\Omega$ and $\beta$ in
Equation \eqref{cond1} as well as on the prequantizing $U(1)$
bundle $Q$ and connection $\sigma$. In general, there is no choice
for which it is injective, as the following example
 shows.
It follows that  the prequantization representation on $H^n_{bas}$
or $H^n_{adm}$ (given by restricting suitably the representation
\eqref{homom}) is generally not faithful for any $n$.

\begin{example}\label{su2}
Consider the Poisson manifold $(S^2\times
\RR^+,\Lambda=t\Lambda_{S^2})$ where $t$ is the coordinate on
$\RR^+$ and $\Lambda_{S^2}$ is the product of the Poisson
structure on $S^2$ corresponding to the standard symplectic form
$\omega_{S^2}$ and the zero Poisson structure on $\RR^+$. (This is
isomorphic to the Lie-Poisson structure on
$\mathfrak{su}(2)^*-\{0\}$.)
 We first
claim that for all choices of $\Omega$ and $A$ in \eqref{preqpois}
(which, as pointed out in Remark \ref{remcond}, is equivalent to
\eqref{cond1}), the $\frac{\partial}{\partial t}$-component of the
vector field $A$ has the form $(ct^2-t)\frac{\partial}{\partial
t}$ for some real constant $c$.

Indeed, notice that $\Lambda+[-t\frac{\partial}{\partial
t},\Lambda]=0$, so
\begin{eqnarray}\label{C}
\tilde{\Lambda}
cp^*\omega_{S^2}=ct^2\Lambda_{S^2}=\Lambda+[A,\Lambda]
\end{eqnarray}
where $A=ct^2\frac{\partial}{\partial t}-t\frac{\partial}{\partial
t}$. Now any vector field $B$ satisfying $[B,\Lambda]=0$ must map
symplectic leaves to symplectic leaves, and since all leaves have
different areas, $B$ must have no $\frac{\partial}{\partial
t}$-component. Hence any vector field satisfying Equation
\eqref{C} has the same $\frac{\partial}{\partial t}$-component as
$A$ above. Now any closed 2-form $\Omega$ on $S^2\times\RR^+ $ is
of the form $cp^*\omega_{S^2}+d\beta$ for some 1-form $\beta$,
where $p:S^2\times \RR^+ \rightarrow S^2$. Since $\tilde{\Lambda}
d\beta=-[\tilde{\Lambda} \beta,\Lambda]$  and $-\tilde{\Lambda}
\beta$ has no $\frac{\partial}{\partial t}$ component, our first
claim is proved.

Now, for any choice of $Q$ and $\sigma$, let $g$ be a function on
$S^2\times \RR^+$ such that $X_{\pi^*g}=X_g^H + (\langle dg,A
\rangle -g)E$ vanishes. This means that  $g$ is a function of $t$
only, satisfying $(ct^2-t)g^{\prime}=g$. For any real number $c$,
there exist non-trivial functions satisfying these conditions, for
example $g=\frac{ct-1}{t}$, therefore for all choices the
homomorphism $g\mapsto -X_{\pi^*g}$ is not injective.

This example also shows that one can not simply omit the vector field
$A$ from the definition of prequantizability, since no choice of $c$
makes $A$ vanish here.
\end{example}

Even though the prequantization representation for functions acting
on  $H^{n}_{adm}$
and  $H^{n}_{bas}$  is usually not faithful for any integer $n$,
 we still have the
following result, which shows that hamiltonian vector fields do act faithfully.
\begin{proposition}\label{inj}
For each integer $n\neq 0$, the map that assigns to an equivalence
class of hamiltonian vector fields $X_{\pi^*g}$ the corresponding
operator on  $H^{n}_{adm}$  or  $H^{n}_{bas}$ is injective.
\end{proposition}
\begin{proof}
Since $H^{n}_{adm}\subset H^{n}_{bas}$, it is enough to consider
the $H^{n}_{adm}$ case.
 Since the hamiltonian vector field of any function on $Q$
is determined up to smooth sections of the singular distribution
$F:=\bar{L}\cap TQ=\{X^H+\langle \alpha, X \rangle E: X\in L\cap
TP\}$, we have to show that, if a $U(1)$-invariant vector field $Y$ on $Q$
annihilates all functions in $H^{n}_{adm}$, then $Y$ must be a
section of $F$.

We start by characterizing the functions in $H^{n}_{adm}$ on
neighborhoods where a constant rank assumption holds:
\begin{lemma}\label{firstlemma}
Let $U$ be an open set in $P$ on which the rank of $L\cap TP$ is
constant and $\bar{U}=\pi^{-1}(U)$. Then a function $\phi$ on
$\bar{U}$ is admissible iff $\phi$ is constant along the leaves of
$F$. Further $\cap_{\phi\in H^{n}_{adm}}\ker d\phi = F$.
\end{lemma}
\begin{proof}We have
\begin{eqnarray}
\phi \text{ admissible } \Leftrightarrow (d\phi,\phi)\subset
\rho_{T^*Q\times \RR}(\bar{L})\Leftrightarrow d\phi \subset
\rho_{T^*Q}(\bar{L}),
\end{eqnarray}
where the first equivalence follows from the formula for
$\bar{L}$, the remark following Definition \ref{dirhamvf} and
 the fact that $\dim (L\cap TP)$ is constant. For any Dirac-Jacobi structure one has
$\rho_{T^*Q}(\bar{L})=(\bar{L}\cap TQ)^{\circ}$, so the first
statement follows.

Now consider the regular foliation of $\bar{U}$ with leaves equal
to $U(1)\cdot \cF$, where $\cF$ ranges over the leaves of
$F|_{\bar{U}}$\footnote{The distribution $F=\bar{L}\cap TQ$ is
clearly involutive; see Definition \ref{defjdbr}.}. Fix $p\in
\bar{U}$ and choose a submanifold $S$ through $p$ which is
transverse to the foliation $U(1)\cdot \cF$. Given any covector
$\xi\in T^*_pS$ we can find a function $\phi$ on $S$ with
differential $\xi$ at $p$, and we extend $\phi$ to $\bar{U}$ so
that it is constant on the leaves of $F$ and equivariant with
respect to the $n$-th power of the standard $U(1)$ action on
$\CC$. Then $\phi$ will lie in $H^{n}_{adm}$ and $d_p \phi$ will
be equal to $\xi$ on $T_pS$,  equal to $2\pi i n$ on $E_p$, and
will vanish on $F_p$. Since we can construct such a function
$\phi\in H^{n}_{adm}$ for any choice of $\xi$, it is clear that a
vector at $p$ annihilated by all functions in $H^{n}_{adm}$ must
lie in $F_p$, so $\cap_{\phi\in H^{n}_{adm}}\ker d\phi \subset F$.
The other inclusion is clear.
\end{proof}

Now we make use of the fact that for any open subset $V$ of $P$
there exists a nonempty open subset $U\subset V$ on which
$\dim(L\cap TP)$ is constant\footnote{Indeed, if $q$ is a point of
$V$ where $\dim(L\cap TP)$ is minimal among all points of $V$, in
a small neighborhood of $q$ $\dim(L\cap TP)$ can not decrease, nor
it can increase because $L\cap TP$ is an intersection of
subbundles.}, and prove Proposition \ref{inj}.

\bigskip
\noindent \textit{End of proof of Proposition \ref{inj}.}~~
Suppose now  the $U(1)$-invariant vector field $Y$ on $Q$ annihilates all functions in
$H^{n}_{adm}$ but is not a section of $F$. Then $Y\notin F$ at all
points of some open set $\bar{U}$. By the remark above, we can
assume that on $\bar{U}$ $\dim(L\cap TP)^H=\dim F$ is constant. By
Lemma \ref{firstlemma} on $\bar{U}$ the vector field $Y$ must be
contained in $F$, a contradiction.
\end{proof}

If we modified the representation \eqref{homom} to act on
\emph{global} admissible or basic functions, the injectivity
statement of Proposition \ref{inj} could fail, as the following
example shows.

\begin{example} Let $P$ be $(\TT^2\times \RR, d\varepsilon)$, where
$\varepsilon=x_3(dx_1+x_3 dx_2)$ with $(x_1,x_2)$ and $x_3$
standard coordinates on the torus and $\RR$ respectively.  This is
a regular presymplectic manifold, so by Lemma \ref{firstlemma} all
basic functions on
any prequantization $Q$ are admissible.
 $P$ is clearly prequantizable, and we can choose
$\Omega=0$ and $\beta=-\rho^*_{TP}\varepsilon$ in the
prequantization condition \eqref{cond1}. Therefore $Q$  is the
trivial $U(1)$ bundle over $P$, with trivial connection
$\sigma=d\theta$ (where $\theta$ is the standard fiber
coordinate). The distribution $F$ on $Q$, as defined at the
beginning of the proof of Proposition \ref{inj}, is one
dimensional, spanned by $2 x_3\frac{\partial}{\partial x_1}
-\frac{\partial}{\partial x_2}-x_3^2\frac{\partial}{\partial
\theta}$. The coefficients $2x_3$, $-1$, and $-x_3^2$ are linearly
independent over $\ZZ$ unless $x_3$ is a quadratic algebraic
integer, so the closures of the leaves of $F$ will be of the form
$\TT^2\times\{x_3\}\times U(1)$ for a dense set of $x_3$'s.
Therefore $C^{\infty}_{adm}(Q,\CC)=C^{\infty}_{bas}(Q,\CC)$
consists exactly of complex functions depending only on $x_3$. For
similar reasons, the admissible functions on $P$ are exactly those
depending only on $x_3$.   But the  vector field $X_{\pi^*g}$ on
$Q$ associated to such a function $g$ has no
$\frac{\partial}{\partial x_3}$
component, so it acts trivially on $C^{\infty}_{adm}(Q,\CC)$.\\
\end{example}

Next we illustrate how the choices involved in the prequantization
representation affect injectivity.

\begin{example} Let $P=S^2\times \RR\times S^1$, with coordinate $t$ on
the $\RR$-factor and $s$ on the $S^1$-factor. Endow $P$ with the
Poisson structure $\Lambda$ which is the product of the zero
Poisson structure on $\RR\times S^1$ and the inverse of an
integral
 symplectic form
$\omega_{S^2}$ on $S^2$. This Poisson manifold is prequantizable;
in Equation \eqref{preqpois} we can choose
$\Omega=p^*\omega_{S^2}$ (where $p:P\rightarrow S^2$) and as $A$
any vector field that preserves the Poisson structure. Each $g \in
C^\infty (P)$ is prequantized by the action of the negative of
1its hamiltonian vector field
$X_{\pi^*g}=(\tilde{\Lambda}dg)^H+(A(g)-g)E$. Therefore the kernel
of the prequantization representation is given by functions of $t$
and $s$ satisfying $A(g)=g$. It is clear that if $A$ is tangent to
the symplectic leaves the representation will be faithful. If $A$
is not tangent to the symplectic leaves, then $A(g)=g$ is an
honest first order differential equation. However, even in this
case the representation might be faithful: it is faithful if we
choose $A=\frac{\partial}{\partial t}$, but not if
$A=\frac{\partial}{\partial s}$.
\end{example}

\begin{remark}\label{variation}
Let $(P,\Lambda)$ be a Poisson manifold such that its symplectic
foliation $\cF$ has constant rank, and assume that $(P,\Lambda)$
is prequantizable (i.e. \eqref{preqpois}, or equivalently
\eqref{cond1}, is satisfied). It follows from the discussion
following Lemma \ref{regular} that, after we fix a prequantizing
$U(1)$-bundle $Q$, the prequantizing Dirac-Jacobi structures on
$Q$ are given by $\bar{L}(Q,\sigma,A)$ where $\sigma$ is fixed and
 $A$ is unique up to vector fields $A'$
preserving $\cF$ such that $\cL_{A'}\Omega_L=0$, i.e. up to vector
fields whose flows are symplectomorphisms between the symplectic
leaves.
 If the topology and geometry of the symplectic leaves of
$P$ ``varies'' sufficiently from one leaf to another (as in
Example \ref{su2} above), then the projection of the $A$'s as
above to $TP/T{\cF}$ will all coincide. Therefore the kernels of
the prequantization representations \eqref{homom}, which associate
to $g \in C^\infty (P)$ the negative of the  hamiltonian vector
field $X_{\pi^*g}=(\tilde{\Lambda}dg)^H+(A(g)-g)E$, will coincide
for all representations arising from prequantizing Dirac-Jacobi
structures over $Q$.
\end{remark}

We end this section with two remarks linked to Kostant's work
\cite{Ko2}.
\begin{remark}\label{hammech}
Kostant (\cite{Ko2}, Theorem 0.1) has observed that the
prequantization of a symplectic manifold can be realized by the
Poisson bracket of a symplectic manifold two dimensions higher,
i.e. that prequantization is ``classical mechanics two dimensions
higher''. In the general context of Dirac manifolds we have seen
in \eqref{homom} that prequantization is given by a Jacobi
bracket\footnote{The bracket on functions on the prequantization
$(Q,\bar{L})$ of a Dirac manifold makes $C^{\infty}_{adm}(Q)$ into
a Jacobi algebra. See Section 5 of \cite{Wa}, which applies
because the constant functions are admissible for the Dirac-Jacobi
structure $\bar{L}$.}; we will now show that Kostant's remark
applies in this context too.

Let $(P,L)$ be a prequantizable Dirac manifold, $(Q,\bar{L})$ its
prequantization and $(Q\times \RR,\tilde{\bar{L}})$ the
``Diracization'' of $(Q,\bar{L})$. To simplify the notation, we
will denote pullbacks of functions (to $Q$ or $Q\times \RR$) under
the obvious projections by the same symbol. Using the homomorphism
\eqref{homJDiracDirac} we can re-write the representation
\eqref{homom} of $C^{\infty}_{adm}(P)$ on
 $C^{\infty}_{adm}(Q,\CC)_{P-loc}$ (or $C^{\infty}_{bas}(Q,\CC)_{P-loc}$) as
$$ g \mapsto e^{-t}\{e^tg,e^t\cdot\}_{Q\times
\RR}=\{e^tg,\cdot\}_{Q\times \RR},$$ i.e. $g$ acts by the Poisson
bracket on $Q\times \RR$.
\end{remark}

\begin{remark}\label{koremark}
Kostant \cite{Ko2} also shows that a prequantizable symplectic
manifold $(P,\Omega)$ can be recovered by reduction from the
symplectization $(Q\times \RR,d(e^t\sigma))$ of its
prequantization $(Q,\sigma)$. More precisely, the inverse of the
natural $U(1)$ action on $Q\times \RR$ is hamiltonian with momentum
map $e^t$, and symplectic reduction at $t=0$ delivers
$(P,\Omega)$. We will show now how to extend this
construction\footnote{Kostant calls the procedure of taking the
  symplectization of the prequantization ``symplectic induction'';
  the term seems to be used here in a different sense from that in
  \cite{KaKS}.} to prequantizable Dirac manifolds.

Let   $(P,L)$, $(Q,\bar{L})$ and   $(Q\times
\RR,\tilde{\bar{L}})$ be as in Remark \ref{hammech}. Since
 $-E\oplus de^t\in
\tilde{\bar{L}}$ we see that $e^t$ is a ``momentum map'' for the
inverse $U(1)$ action  on $Q\times \RR$, and by Dirac reduction
\cite{BlR} at the regular value 1 we obtain $L$: indeed,
 the
pullback of $\tilde{\bar{L}}$ to $Q\times \{0\}$ is easily seen to
be $\{(X^H+(\langle X\oplus\xi,\beta \rangle -g)E)\oplus \pi^*\xi:
X\oplus \xi\in L\}$, and its pushforward via $\pi:Q\rightarrow P$
is exactly $L$.
\end{remark}

\section{The line bundle approach}\label{line}

In this section we will prequantize a Dirac manifold $P$ by
letting its admissible functions act on sections of a hermitian
line bundle
$K$ over $P$. This approach was first taken by Kostant for
symplectic manifolds and was extended by Huebschmann \cite{Hu} and
Vaisman \cite{Va2}
  to Poisson manifolds. The construction of this section
generalizes Vaisman's and turns out to be  equivalent to the one
we described in Sections \ref{spaces} and \ref{representation}.

\begin{definition} \cite{F} \label{Lconn}
Let $(A,[\cdot,\cdot],\rho)$ be a Lie algebroid over the manifold
$M$ and $K$ a real vector bundle over $M$. An $A$-connection on
the vector bundle $K\rightarrow M$ is a map $D: \Gamma(A) \times
\Gamma(K) \rightarrow \Gamma(K)$ which is $C^{\infty}(M)$-linear
in the $\Gamma(A)$ component and satisfies
$$D_e(h\cdot s)=h\cdot D_es+\rho e (h)\cdot s,$$
for all $e\in \Gamma(A)$, $s\in \Gamma(K)$ and $h\in
C^{\infty}(M)$. The curvature of the $A$-connection is the map
$\Lambda^2 A^* \rightarrow End(K)$ given by
$$R_D(e_1,e_2)s=D_{e_1}D_{e_2} s -D_{e_2}D_{e_1}s
-D_{[e_1,e_2] }s .$$

If $K$ is a complex vector bundle, we define an $A$-connection on
$K$ as above, but with $C^{\infty}(M)$ extended to the
complex-valued smooth functions.

\end{definition}
\begin{remark} When $A=TM$ the definitions above specialize to the
usual notions of covariant derivative and curvature. Moreover,
given an ordinary connection $\nabla$ on $K$, we can pull it back
to a $A$-connection by setting $D_e=\nabla_{\rho e}$.\end{remark}

With this definition we can easily adapt Vaisman's construction
\cite{Va2} \cite{Va3}, extending it  from the case where $L=T^*P$
is the Lie algebroid of a Poisson manifold to the case where $L$
is a Dirac structure. We will act on locally defined, basic
sections.

\begin{lemma}\label{line1}
 Let $(P,L)$ be
a Dirac manifold and $K$ a hermitian line bundle over $P$ endowed
with an $L$-connection
 $D$. Then $R_D=2\pi i \Upsilon$, where
$\Upsilon =\langle \cdot,\cdot \rangle_-|_L$, iff the
correspondence
$$\hat{g}s=-(D_{X_g\oplus dg}s+2\pi igs)$$ defines a Lie algebra
representation of $C^{\infty}_{adm}(P)$ on $\{s\in\Gamma(K)_{loc}:
D_{Y\oplus 0}s=0 \text{ for } Y\in L\cap TP\}$, where $X_g$ is any
choice of hamiltonian vector field for $g$.
\end{lemma}
\begin{proof}
If $\hat{g}$ and $s$ are as above, then clearly $\hat{g}s$ is a
well-defined section of $K$. We will now show that $\hat{g}s \in
\{s\in\Gamma(K)_{loc}: D_{Y\oplus 0}s=0 \text{ for }Y\in L\cap
TP\}$, so that the above ``representation'' is well-defined. The
case where $Y\in L\cap T_pP$ can be locally extended to a smooth
section of $L\cap TP$ is easy, whereas the techniques (see Section
2.5 of \cite{F}) needed for general case are much more involved.

 The section $X_g\oplus dg$ of $L$ induces a flow $\phi_t$ on $P$
(which is just the flow of the vector field $X_g$) and a
one-parameter family of bundle automorphisms $\Phi_t$ on $TP\oplus
T^*P$ which (see Section 2.4 in \cite{Co}) preserves $L$, and
which takes $L$-paths to $L$-paths\footnote{For any algebroid $A$
over $P$ an $A$-path is a defined as a path $\Gamma(t)$ in $A$
such that the anchor maps $\Gamma(t)$ to the velocity of the base
path $\pi(\Gamma(t))$.}. Further, $\Phi_t$ acts on the sections
$s$ of the line bundle $K$ too, as follows: $(\Phi_t^*s)_p$ is the
parallel translation of $s_{\phi_t(p)}$ along the $L$-path
$\Phi_{\bullet}(X_g\oplus dg)_p=(X_g\oplus
dg)_{\phi_{\bullet}(p)}$. Now $(D_{(X_g\oplus
dg)}s)_p=\frac{\partial}{\partial t}|_0 (\Phi_t^*s)_p$, and
$(D_{(Y\oplus 0)}D_{X_g\oplus dg}s)_p=\frac{\partial}{\partial
t}|_0 (D_{(Y\oplus0)}\Phi_t^*s)_p.$ For every $t$, since $\phi_t$
preserves $L\cap TP$, we have
\begin{eqnarray}\label{zeroinK}0=(D_{({\phi_t}_*Y\oplus 0)}s)_{\phi_t(p)}=
\frac{\partial}{\partial \epsilon}\Big|_0 \para^{\epsilon}_0
s_{\phi_t(\gamma(\epsilon))}\end{eqnarray}
 where $\Gamma$ is an
$L$-path starting at $(Y\oplus 0)\in L_p$, $\gamma$ is its base
path, and $\para^{\epsilon}_0$  is parallel translation along the
$L$-path $\Phi_t(\Gamma(\bullet))$. (This notation denotes the
path $\epsilon \mapsto \Phi_t(\Gamma(\epsilon))$.)
Now we parallel translate the
element \eqref{zeroinK} of $K_{\phi_t(p)}$ to $p$ using the
$L$-path $\Phi_{\bullet}(X_g\oplus dg)_p$, and compare the result
with
\begin{eqnarray}\label{otherterm}(D_{(Y\oplus0)}\Phi_t^*s)_p=\frac{\partial}{\partial
\epsilon}\Big|_0
\para^{\epsilon}_0 \para^{t}_0
s_{\phi_t(\gamma(\epsilon))},\end{eqnarray}
 where
the parallel translation is taken first along
$\Phi_{\bullet}(X_g\oplus dg)_{\gamma(\epsilon)}$ and then along
$\Gamma(\bullet)$.

The difference  between \eqref{otherterm} and the parallel
translation to $p$ of \eqref{zeroinK} lies only in the order in
which the parallel translations are taken. Now applying
$\frac{\partial}{\partial t}|_0$ to this difference (and recalling
that $\Phi_t(X_g\oplus dg)_p=(X_g\oplus dg)_{\phi_t(p)}$) we
obtain the evaluation at $p$ of
$$D_{\Phi_t\Gamma(\epsilon)}D_{(X_g\oplus dg)}s- D_{(X_g\oplus
dg)} D_{\Phi_t\Gamma(\epsilon)}s,$$ which by the definition of
curvature is just $$(D_{[\Phi_t\Gamma(\epsilon),X_g\oplus
dg]}s)_p+\Upsilon(Y\oplus0,(X_g\oplus dg)_p)s.$$ The second term
vanishes because $Y\in L\cap T_pP$, and using the fact that
$\Phi_t$ is the flow generated by $X_g$ one sees that the Courant
bracket in the first term is also zero. Altogether we have proven
that $(D_{(Y\oplus 0)}D_{X_g\oplus dg}s)_p$ vanishes, and from
this is follows easily that the ``representation'' in the
statement of the lemma is well defined.

Since
$$[\hat{f},\hat{g}]=D_{X_f\oplus df}D_{X_g\oplus dg}-D_{X_g\oplus dg}D_{X_f\oplus df}+2
\pi i(X_f(g)-X_g(f)),$$ using $-[X_f\oplus df,X_g\oplus
dg]=X_{\{f,g\}} \oplus d\{f,g\}$
 (\cite{Co}, Prop. 2.5.3) we see that the condition on $R_D$
 holds iff $[\hat{f},\hat{g}]=\widehat{\{f,g\}}$.
\end{proof}
Now assume that the prequantization condition
\eqref{cond1} is satisfied, i.e. that there exists a closed
integral 2-form $\Omega$ and a Lie algebroid 1-cochain $\beta$ for
such that
$$\rho_{TP}^*\Omega=\Upsilon+d_L \beta.$$ Then
 we can construct  an $L$-connection $D$ satisfying the
property of the previous lemma:
\begin{lemma}\label{line2}
Let $(A,[\cdot,\cdot],\rho)$ be a Lie algebroid over the manifold
$M$, $\Omega$ a closed integral 2-form on $M$, and $\nabla$ a
connection (in the usual sense)
on a hermitian line bundle $K$ with curvature
$R_{\nabla}=2\pi i \Omega$. If $\rho^*\Omega=\Upsilon+d_L \beta$
for a 2-cocycle $\Upsilon$ and a 1-cochain $\beta$ on $A$, then
the $A$-connection $D$ defined by
$$D_e=\nabla_{\rho e}-2 \pi i \langle e,\beta \rangle$$
has curvature $R_D=2\pi i \Upsilon$.
\end{lemma}
\begin{proof}
An easy computation shows
$$R_D(e_1,e_2)=R_{\nabla}(\rho e_1,\rho e_2) +2 \pi i
(-\rho e_1 \langle e_2,\beta \rangle+\rho e_2 \langle e_1,\beta
\rangle+ \langle [e_1,e_2],\beta \rangle),$$ which using
$\rho^*\Omega=\Upsilon+d_L \beta$ reduces to $2\pi i
\Upsilon(e_1,e_2)$.
\end{proof}

Altogether we obtain that
\begin{equation*}\label{lbhomom}\hat{g}=-[\nabla_{X_g}-2
\pi i(\langle X_g\oplus dg,\beta \rangle-g)] \end{equation*}
determines a representation of $C^{\infty}_{adm}(P)$ on
$\{s\in\Gamma(K)_{loc}: \nabla_Y s-2 \pi i \langle Y\oplus 0,\beta
\rangle s=0 \text{ for }Y\in L \cap TP\}$. Notice that, when $P$
is symplectic, we recover Kostant's prequantization mentioned in
the introduction.
 Now let
$Q\rightarrow P$ be the $U(1)$-bundle corresponding to $K$, with
the connection form $\sigma$ corresponding to $\nabla$. If
$\bar{s}$ is the $U(1)$-antiequivariant complex valued function on
$Q$ corresponding to the section $s$ of $K$, then $X^H(\bar{s})$
corresponds to $\nabla_Xs$ and $E(\bar{s})$ to $-2\pi i s$. Here
$X\in TP$, $X^H \in \ker \sigma$ its horizontal lift to $Q$, and
$E$ is the infinitesimal generator of the $U(1)$ action on $Q$ (so
$\sigma(E)=1$). Translating the above representation to the $U(1)$-bundle picture, we see that $\hat{g}=-[X_g^H+(\langle X_g\oplus
dg,\beta \rangle-g)E]$ defines a representation of
$C^{\infty}_{adm}(P)$ on
\begin{equation*}\begin{split}
\{&\bar{s}\in C^{\infty}(Q, \CC)_{P-loc}: \bar{s} \text{ is
}U(1)\text{-antiequivariant and }\\(&Y^H+\langle Y\oplus 0, \beta
\rangle E)\bar{s}=0 \text{ for }Y\in L\cap
TP\},\end{split}\end{equation*} which is nothing else than
$H^{-1}_{bas}$ as defined in Section \ref{representation}.
 Since $X_g^H+(\langle X_g\oplus dg,\beta \rangle-g)E$ is the
hamiltonian vector field of $\pi^*g$ (with respect to the
Dirac-Jacobi structure $\bar{L}$ on $Q$ as in Theorem
\ref{thmpreq}), we see that this is exactly
 our prequantization representation given
by Equation \eqref{homom} restricted to $H^{-1}_{bas}$.

\subsection{Dependence of the prequantization on choices: the line
bundle point of view}

In Subsection \ref{choices} we gave a classification the
Dirac-Jacobi structures induced on the prequantization of a given
Dirac manifold, and hence also a classification of the
corresponding prequantization representations. Now we will see
that the line bundle point of view allows for an equivalent but
clearer classification.

Recall that, given a Dirac manifold satisfying the prequantization
condition \eqref{cond1}, we associated to it a hermitian
line bundle $K$ and
a representation as in Lemma \ref{line1}, where the $L$-connection
$D$  is given as in Lemma \ref{line2}

\begin{proposition}
Fix a line bundle $K$ over $P$ with $(\rho_{TP}^* \circ
i_*)c_1(K)=[\Upsilon]$. Then all the hermitian  $L$-connections of $K$ with
curvature $\Upsilon$ are given by the $L$-connections constructed
in Lemma \ref{line2}. Therefore there is a surjective map from
the set of isomorphism classes of prequantization
representations of $(P,L)$ to the space
$(\rho_{TP}^* \circ i_*)^{-1}[\Upsilon]$ of topological types; the set with a given type is a principal homogeneous space for $H_L^1(P,U(1)).$
\end{proposition}
\begin{proof}
 Exactly as in the case of ordinary connections one shows that
the difference of two hermitian  $L$-connections on $K$ is a
section of $L^*$, whose $d_L$-derivative is the difference of the
curvatures. Fix a choice of $L$-connection $D$ as in Lemma
\ref{line2}, say given by $D_{(X\oplus\xi)}=\nabla_{X}-2 \pi i
\langle X\oplus \xi ,\beta \rangle$. Another $L$-connection $D'$
with curvature $\Upsilon$ is given by
$D'_{(X\oplus\xi)}=\nabla_{X}-2 \pi i \langle X\oplus \xi
,\beta+\beta' \rangle$ for some $d_L$-closed section $\beta '$ of
$L^*$, hence it arises as in Lemma \ref{line2}. This shows the
first claim of the proposition. Since, as we have just seen, the
$L$-connections with given curvature differ by $d_L$-closed
sections of $L^*$ and since $U(1)$-exact sections of $L^*$ give
rise to gauge equivalences of hermitian line bundles with
connections, the second claim follows as well.
\end{proof}

Using Lemma \ref{lemmathreeNEW} it is easy to see that choices of
$(\sigma,\beta)$
giving rise to the same $L$-connection (as in Lemma \ref{line2})
also give rise to the same Dirac-Jacobi structure $\bar{L}$, in
accord with the results of Section \ref{choices}.
 Given this, it is  natural to try to express the
Dirac-Jacobi structure $\bar{L}$ intrinsically in terms of the
$L$-connection to which it corresponds; this is subject of work in
progress.

\section{Prequantization of
Poisson and Dirac structures associated to  contact manifolds}
\label{seclebrun} We have already mentioned in Remark
\ref{koremark} the symplectization construction, which associates
to a manifold $M$ with contact form $\sigma$ the manifold
$M\times\RR$ with symplectic form $d(e^t\sigma)$.  The
construction may also be expressed purely in terms of the
cooriented contact distribution $C$ annihilated by $\sigma$.  In
fact, given any contact distribution, its nonzero annihilator
$C^\circ$ is a (locally closed) symplectic submanifold of $T^*M$.
When $C$ is cooriented, we can select the positive component
$C^\circ_+$.  Either of these symplectic manifolds is sometimes
known as the symplectization of $(M,C)$.  It is a bundle over $M$
for which a trivialization (which exists in the cooriented case)
corresponds to the choice of a contact form $\sigma$ and gives a
symplectomorphism between this ``intrinsic'' symplectization and
$(M\times\RR,d(e^t \sigma))$. The contact structure on $M$ may be
recovered from its symplectization along with the conformally
symplectic $\RR$ action generated by $\partial/\partial t$.

One may partially compactify $C^\circ_+$ (we stick to the
cooriented case for simplicity) at either end to get a manifold
with boundary diffeomorphic to $M$.  The first, and simplest way,
is simply to take its closure $C^\circ_{0,+}$ in the cotangent
bundle by adjoining the zero section.  The result is a
presymplectic manifold with boundary, diffeomorphic to
$M\times[0,\infty)$ with the exact 2-form $d(s\sigma)=ds\wedge
\sigma + s d\sigma$, where $s$ is the exponential of the coordinate
$t$ in $\RR$.
For positive $s$, this is
symplectic; the characteristic distribution of $C^\circ_{0,+}$
 lives along the boundary $M\times \{0\}$, where it may
be identified with the contact distribution $C$.  This is
highly nonintegrable even though $d(s\sigma)$ is closed,
so we have another example of the phenomenon alluded to in the
discussion after Definition \ref{defdirstr}.

We also note that the basic functions
on $C^\circ_{0,+}$ are just those which are constant on
$M\times \{0\}$.  One can prove that all of these functions are
admissible as well, even though the characteristic distribution is
singular.  It would be interesting to characterize the Dirac
structures for which these two classes of functions coincide.

To compactify the other end of  $C^\circ_{+}$, we begin by identifying $C^\circ_+$ with
the positive part of its dual $(TM/C)_+$, using the ``inversion''
map $j$ which takes $\phi \in C^\circ_+$ to the unique element
$X\in (TM/C)_+$ for which $\phi(X)=1$.  We then form the union
$C^\circ_{+,\infty}$ of $C^\circ_+$ with the zero section in
$TM/C$ and give it the topology and differentiable structure
induced via $j$ from the closure of $(TM/C)_+$.  It was discovered
by LeBrun \cite{Le} that the Poisson structure on $C^\circ_+$
corresponding to its symplectic structure extends smoothly to
$C^\circ_{+,\infty}$.  We call $C^\circ_{+,\infty}$ with this
Poisson structure the {\bf
  LeBrun-Poisson manifold} corresponding to the contact manifold
$(M,C)$.

To analyze the LeBrun-Poisson structure more closely, we introduce the
inverted coordinate $r=1/s$, which takes values in $[0,\infty)$ on
  $C^\circ_{+,\infty}$.  In suitable local coordinates on $M$, the
  contact form $\sigma$ may be written as $du+\sum p_i dq^i$.  On the
  symplectization, we have the form $d(r^{-1}(du+\sum p_i dq^i))$.
 The
  corresponding Poisson structure turns out to be
$$\Lambda=r\left[\left(r\frac{\partial}{\partial
r} + \sum p_{i}\frac{\partial}{\partial p_{i}}\right) \wedge
\frac{\partial}{\partial u} + \sum \frac{\partial}{\partial q^{i}}
\wedge \frac{\partial}{\partial p_{i}}\right].$$
 From this formula we  see not only that $\Lambda$ is smooth at
$r=0$ but also that its linearization
$$r \sum \frac{\partial}{\partial q^{i}}
\wedge \frac{\partial}{\partial p_{i}}$$ at the origin (which is a
``typical'' point, since $M$ looks the same everywhere) encodes the
contact subspace in terms of the symplectic leaves in the tangent Poisson
structure.

We may take the union of the two compactifications above to get a
manifold $C^\circ_{0,+,\infty}$ diffeomorphic to $M$ times a
closed interval.  It is presymplectic at the $0$ end and Poisson
at the $\infty$ end, so it can be treated globally only as a Dirac
manifold. In what follows, we will simply denote this Dirac
manifold as $(P,L)$.

To prequantize $(P,L)$, we first notice that its Dirac structure
is ``exact'' in the sense that the cohomology class $[\Upsilon]$
occurring in the condition \eqref{cond0} is zero.  In fact, on the
presymplectic end, $L$ is isomorphic to $TP$, and $\Upsilon$ is
identified with the form $d(s\sigma)$, so we can take the cochain
$\beta$ to be the section of $L^*$ which is identified with
$-s\sigma$.  To pass to the other end, we compute the projection
of this section of $L^*$  into $TP$ and find that it is just the Euler vector field
$A=s\frac{\partial}{\partial s}$. In terms of the inverse
  coordinate $r$, $A = -r\frac{\partial}{\partial r}$.
  (The reader may check that the Poisson differential of this vector field is
  $-\Lambda$, either by direct computation or using
the degree $1$ homogeneity of $\Lambda$ with respect to  $r$.) On
the Poisson end, $L^*$ is isomorphic to $TP$,
so $-r\frac{\partial}{\partial r}$
defines a smooth continuation of $\beta$ to all of $P$.

Continuing with the prequantization, we can take the 2-form
$\Omega$ to be zero and the $U(1)$-bundle $Q$ to be the product
$P\times U(1)$ with the trivial connection $d\theta$, where
$\theta$ is the ($2\pi$-periodic) coordinate on $U(1)$. On the
presymplectic end, the Dirac-Jacobi structure is defined by the
1-form $\sigma =  s\sigma + \theta$, which is a contact form when
$s\neq 0$.

  On the Poisson
end, we get the Jacobi structure $(\Lambda^H + E\wedge A^H, E)$ which
in coordinates becomes
\begin{equation}
\label{lebrunjacobi}
\left(r\left[\left(r\frac{\partial}{\partial
r} + \sum p_{i}\frac{\partial}{\partial p_{i}}\right) \wedge
\frac{\partial}{\partial u} + \sum \frac{\partial}{\partial q^{i}}
\wedge \frac{\partial}{\partial
  p_{i}}-\frac{\partial}{\partial\theta}\wedge
\frac{\partial}{\partial r} \right] ,\frac{\partial}{\partial\theta}
 \right).
\end{equation}

\section{Prequantization by circle actions with fixed points}\label{fixed}
Inspired by a construction of Engli\v{s} \cite{en:weighted}
in the complex
setting, we modify the prequantization in the previous section
 by ``pinching''
the boundary component $M\times U(1)$ at the Poisson end and replacing it
by a copy of $M$.  To do this, we identify $U(1)$ with the unit
circle in the plane $\RR^2$ with coordinates $(x,y)$.  In
addition, we make a choice of contact form on $M$ so that $P$ is
identified with $M\times [0,\infty]$, with the coordinate $r$ on
the second factor.  Next we choose a smooth nonnegative real
valued function $f:[0,\infty]\to\RR$ such that, for some $\epsilon
> 0$, $f(r)=r$ on $[0,\epsilon]$ and $f(r)$ is constant on
$[2\epsilon,\infty]$. Let $Q'$ be the submanifold of
$P\times\RR^2$ defined by the equation $x^2+y^2=f(r)$.

Radial projection in the $(x,y)$ plane determines a map $F:Q\to
Q'$ which is smooth, and in fact a diffeomorphism, where $r>0$.
The boundary $M\times U(1)$ of $Q$ is projected smoothly to
$M\times (0,0)$ in $Q'$, but $F$ itself is not smooth along the
boundary.  We may still use $F$ to transport the Jacobi structure
on $Q$ to the part of $Q'$ where $r>0$.  For small $r$, we have
$x=\sqrt{r} \cos \theta$ and $y=\sqrt{r} \sin\theta$, so
$r=x^2+y^2$, $r\frac{\partial}{\partial r}=
\frac{1}{2}(x\frac{\partial}{\partial x}+y\frac{\partial}{\partial
y}),$ and
$\frac{\partial}{\partial\theta}=x\frac{\partial}{\partial y} -
y\frac{\partial}{\partial x} .$ Using these substitutions to write
the Jacobi structure \eqref{lebrunjacobi} with polar coordinates
$(r,\theta)$ replaced by rectangular coordinates $(x,y)$, we see
immediately that the structure extends smoothly to a Jacobi
structure on the Poisson end of $Q'$ and to a Dirac-Jacobi
structure on all of $Q'$, and that the projection $Q'\to P$, like
$Q\to P$ pushes the Dirac-Jacobi structure on $Q'$ to the Dirac
structure on $P$.   (Thus, the projection is a ``forward
Dirac-Jacobi map''; see the beginning of Section \ref{secjacdir}.)
The essential new feature of $Q'$
 is that the vector field
$E'=x\frac{\partial}{\partial y}-y\frac{\partial}{\partial x}$ of the
Jacobi structure on $Q'$ vanishes along the locus $x=y=0$ where the
projection is singular.

The vanishing of $E'$ at some points means that the Jacobi
structure on $Q'$ does not arise from a contact form, even on the
Poisson end, where $r < \infty$.  However, it turns out that we
can turn it into a contact structure by making a conformal change,
i.e. by multiplying the bivector by $1/f$ and replacing $E'$ by
$E'/f + X_{1/f}$.  The resulting Jacobi structure still extends
smoothly over $Q'$, and now comes from a contact structure over
the Poisson end; the price we pay is that the projection to $P$ is
now a conformal Jacobi map rather than a Jacobi map.

\begin{remark}
Looking back at the construction above, we see that we have embedded
any given contact manifold $M$ as a codimension 2 submanifold in
another contact manifold.  Our construction depended only on the choice of
a contact form.  On the other hand, Eliashberg and Polterovich
\cite{el-po:partially} construct a similar embedding in a canonical way,
without the choice of a
contact form.   It is not hard to show that the choice of a contact
form defines a canonical isomorphism between our contact manifold and
theirs.
\end{remark}

\begin{example}
Let  $M$ be the unit sphere in $\CC^n$, with the contact structure
induced from the Cauchy-Riemann structure on the boundary of the
disc $D^{2n}$.  It turns out that a neighborhood $U$ of $M$ in the disc
can be mapped diffeomorphically to a neighborhood $V$ of $M$ at the
Poisson end in its LeBrun-Poisson manifold $P$ so that the symplectic
structure on the interior of $V$ pulls back to the
symplectic structure on $U$ coming from
the K\"ahler structure on the open disc, viewed as complex hyperbolic space.
  If we now pinch the end of the prequantization $Q$, as above, the
part of the contact manifold $Q'$ lying over $V$ can be glued to the
usual prequantization of the open disc so as to obtain a compact
contact manifold $Q''$ projecting by a ``conformal Jacobi map'' to the
closed disc.  The fibres of the map are the orbits of a $U(1)$-action
which is principal over the open disc.  In fact,
$Q''$ is just the unit sphere in $\CC^{n+1}$
with {\em its} usual contact structure.
 All this is the symplectic analogue of the complex construction
by Engli\v{s} \cite{en:weighted}, who enlarges a bounded
pseudoconvex domain $D$ in $\CC^n$ to one in $\CC^{n+1}$ with a
$U(1)$ action on its boundary which degenerates just over the
boundary of $D$.
\end{example}

The ``moral'' of the story in this section is that, in prequantizing a
Poisson manifold $P$ whose Poisson structure degenerates along a submanifold, one might want to allow the prequantization
bundle to be a Jacobi manifold $Q$ whose vector field $E$ generates a
$U(1)$ action having fixed points and for which the quotient
projection $Q\to P$ is a Jacobi map.

\section{Final remarks and questions}\label{remarks}
We conclude with some suggestions for further research along the lines
initiated in this paper.

\subsection{Cohomological prequantization}
Cohomological methods have already been used in geometric
quantization of symplectic manifolds: rather than the space of
global polarized sections, which may be too small or may have
other undesirable properties, one looks at the higher cohomology
of the sheaf of local polarized sections.  (An early reference on
this approach is \cite{sn:cohomology}.)  When we deal with Dirac
(e.g. presymplectic) manifolds, it may already be interesting to
introduce cohomology at the prequantization stage.  There are two
ways in which this might be done.

The first approach, paralleling that which is done with polarizations,
is to replace the Lie algebra of global admissible functions on a
Dirac manifold $P$ by the cohomology of the sheaf of Lie algebras
of local admissible functions.  Similarly, one would replace the
sheaf  of $P$-local functions on $Q$ by its cohomology.
The first sheaf cohomology
 should then act on the second.

The other approach, used by Cattaneo and Felder \cite{CF} for the deformation
quantization of coisotropic submanifolds of Poisson manifolds,
would apply to Dirac manifolds $P$ whose characteristic
distribution is regular. Here, one introduces the ``longitudinal
de Rham complex'' of differential forms along the leaves of the
characteristic foliation on $P$. The zeroth cohomology of this
foliation is just the admissible functions, so it is natural to
consider the full cohomology, or even the complex itself.  It
turns out that, if one chooses a transverse distribution to the
characteristic distribution, the transverse Poisson structure
induces the structure of an $L_\infty$ algebra on the longitudinal
de Rham complex.  Carrying out a similar construction on a
prequantization $Q$ should result in an $L_\infty$ representation
of this algebra.

\subsection{Noncommutative prequantization}
If the characteristic distribution of a Dirac structure $P$ is
regular, we may consider the groupoid algebra associated to the
characteristic foliation as a substitute for the admissible functions.
By adding some extra structure, as in \cite{BG}\cite{Ta}\cite{Xu},
 we can make this groupoid
algebra into a noncommutative Poisson algebra.  This means that the
Poisson bracket is not a Lie algebra structure, but rather a class
with degree 2 and square 0 in the Hochschild cohomology of the
groupoid algebra.  It should be interesting to define a notion of
representation for an algebra with such a cohomology class, and to
construct such representations from prequantization spaces.  Such a
construction should be related to the algebraic quantization of Dirac
manifolds introduced in \cite{TW}.



\noindent Alan Weinstein\\
      Dept. of Mathematics, University of California\\
      Berkeley, CA  94720, U.S.A.\\
  alanw@math.berkeley.edu

\vspace*{0.8cm}

\noindent Marco Zambon\\
      Mathematisches Institut, Universit\"at Z\"urich\\
      Winterthurerstr. 190, 8057 Z\"urich, Switzerland\\
  zambon@math.unizh.ch

\label{lastpage}
\end{document}